\documentclass[11pt]{article} 

\usepackage{iclr2026_conference,times}
\usepackage[pagewise]{lineno}

\usepackage{amsmath,amsfonts,bm}

\def\eqref#1{equation~\ref{#1}}

\def\1{\bm{1}}

\def\rb{{\textnormal{b}}}

\def\rf{{\textnormal{f}}}

\def\rv{{\textnormal{v}}}

\def\rx{{\textnormal{x}}}
\def\ry{{\textnormal{y}}}

\def\rvu{{\mathbf{i}}}

\def\rvu{{\mathbf{u}}}
\def\rvv{{\mathbf{v}}}
\def\rvw{{\mathbf{w}}}
\def\rvW{{\mathbf{W}}}
\def\rvx{{\mathbf{x}}}
\def\rvy{{\mathbf{y}}}

\def\rmX{{\mathbf{X}}}

\def\vzero{{\bm{0}}}

\def\va{{\bm{a}}}
\def\vb{{\bm{b}}}

\def\vf{{\bm{f}}}
\def\vg{{\bm{g}}}

\def\vr{{\bm{r}}}

\def\vu{{\bm{u}}}

\def\vx{{\bm{x}}}
\def\vy{{\bm{y}}}

\def\mA{{\bm{A}}}
\def\mB{{\bm{B}}}

\def\mD{{\bm{D}}}

\def\mI{{\bm{I}}}

\def\mQ{{\bm{Q}}}

\def\mLambda{{\bm{\Lambda}}}

\DeclareMathAlphabet{\mathsfit}{\encodingdefault}{\sfdefault}{m}{sl}
\SetMathAlphabet{\mathsfit}{bold}{\encodingdefault}{\sfdefault}{bx}{n}

\def\gH{{\mathcal{H}}}

\def\sF{{\mathbb{F}}}

\def\sH{{\mathbb{H}}}

\newcommand{\E}{\mathbb{E}}

\newcommand{\R}{\mathbb{R}}

\newcommand{\Var}{\mathrm{Var}}

\DeclareMathOperator*{\argmin}{arg\,min}

\usepackage{hyperref}       
\usepackage{url}            
\usepackage{amssymb, amsthm, mathtools}
\usepackage{graphicx}
\usepackage{subcaption}
\usepackage{hyperref}
\usepackage[utf8]{inputenc} 
\usepackage[T1]{fontenc}    
\usepackage{booktabs}       
\usepackage{nicefrac}       
\usepackage{microtype}      
\usepackage{xcolor}         
\usepackage{float}          
\usepackage[normalem]{ulem}
\usepackage{wrapfig}        
\DeclareMathOperator{\dif}{d \!}        
\newcommand{\dt}{\dif t}
\newcommand{\zero}{\mathbf{0}}
\newcommand{\balpha}{\bm\alpha}
\newcommand{\psivec}{\bm\psi}
\newcommand{\phivec}{\bm\phi}
\newcommand{\alphavec}{\bm\alpha}
\def\rvU{{\mathbf{U}}}
\DeclareMathOperator{\spn}{span}
\newcommand{\inp}[1]{\langle #1 \rangle}
\newcommand{\norm}[1]{||#1||}
\newcommand{\loss}{\mathcal{E}}
\newcommand{\normal}{\mathcal{N}}
\DeclareMathOperator{\Dim}{dim}

\newtheorem{theorem}{Theorem}

\newtheorem{lemma}{Lemma}
\newtheorem{definition}{Definition}
\title{Noise-Aware System Identification for High-Dimensional Stochastic Dynamics}

\author{Ziheng Guo \\
Department of Mathematics\\
University of Houston\\
Houston, TX $77204$, USA \\
\texttt{zguo24@cougarnet.uh.edu} \\
\And
Igor Cialenco \\
Department of Applied Mathematics \\
Illinois Institute of Technology \\
Chicago, IL $60616$, USA\\
\texttt{cialenco@illinoistech.edu} \\
\And
Ming Zhong \\
Department of Applied Mathematics \\
University of Houston \\
Houston, TX $77204$, USA\\
\texttt{mzhong3@central.uh.edu} \\
}

\iclrfinalcopy 
\begin{document}

\maketitle

\vspace{-2em}
\begin{center}
First circulated: October $10$, $2024$ \\
This version: March $8$, $2026$ 
\end{center}

\begin{abstract}
Stochastic dynamical systems are ubiquitous in physics, biology, and engineering, where both deterministic drifts and random fluctuations govern system behavior. Learning these dynamics from data is particularly challenging in high-dimensional settings with complex, correlated, or state-dependent noise. We introduce a noise-aware system identification framework that jointly recovers the deterministic drift and full noise structure directly from the trajectory data, without requiring prior assumptions on the noise model. Our method accommodates a broad class of stochastic dynamics, including colored and multiplicative noise, that scales efficiently to high-dimensional systems, and accurately reconstructs the underlying dynamics. Numerical experiments on diverse systems validate the approach and highlight its potential for data-driven modeling in complex stochastic environments.
\end{abstract}

\section{Introduction}\label{sec:intro}
Stochastic differential equations (SDEs) provide a fundamental and versatile framework for modeling systems in which random fluctuations are intrinsic to the dynamics~\citep{Evans, Sarkka_Solin_2019}. Compared to deterministic ordinary differential equations (ODEs), SDEs incorporate noise explicitly--often through a Brownian motion term--allowing them to capture variability and uncertainty that strongly influence system behavior. This capability is essential for representing complex phenomena in physics, biology, chemistry, and finance, where stochasticity can be a dominant factor. By incorporating deterministic forces and random fluctuations in a unified mathematical description, SDEs offer a flexible modeling approach that is both theoretically rigorous and practically relevant.

We consider SDEs of the form
\[
d\rvx_t = \vf(\rvx_t)\dif t + \sigma(\rvx_t) \dif \rvw_t, \quad \rvx_t, \rvw_t \in \R^D,
\]
where the drift $\vf: \R^D \rightarrow \R^D$ and the diffusion coefficient $\sigma: \R^D \rightarrow \R^{D \times D}$ are potentially unknown.  The driving noise $\rvw_t$ is a vector of independent standard Brownian motions. The noise structure of the SDE system is described by a state dependent covariance matrix $\Sigma: \R^D \rightarrow \R^{D \times D}$, where $\Sigma = \sigma \sigma^\intercal$.  This general formulation encompasses many classical and modern models. In physics, the Langevin equation~\citep{sachs2017langevin, doi:10.1142/8195, EGNS2008, Talay2002StochasticHS} describes microscopic particle dynamics under both systematic forces and thermal fluctuations. In biology, stochastic Lotka–Volterra models~\citep{takeuchi2006evolution} capture population interactions in fluctuating environments, while other SDE-based models describe cellular processes and gene expression noise~\citep{SZEKELY201414, DP2011}. In chemistry, the chemical Langevin equation~\citep{10.1063/1.4948407} accounts for reaction kinetics in small-molecule regimes, where random molecular collisions cannot be neglected. In finance, SDEs form the basis of models such as Black–Scholes~\citep{0b9b8115-a8b8-3422-8e1c-a62077de6621, hull2017options}, Vasicek~\citep{VASICEK1977177}, and Heston~\citep{10.1093/rfs/6.2.327}, which incorporate uncertainty in asset prices, interest rates, and volatility. More recently, SDE formulations have emerged as the mathematical backbone of diffusion models in machine learning~\citep{ho2020denoising, song2020score}, enabling state-of-the-art generative modeling methods.

Accurate application of SDEs requires careful calibration to empirical data so that both the deterministic drift and stochastic noise are faithfully represented. This is crucial for predictive power and for preserving physical interpretability. In many traditional settings, the functional forms of $\vf$ and $\sigma$ are assumed known up to a small set of parameters, which can be estimated via least squares or related regression techniques~\citep{MrazekPospisil2017, 935092}. However, in modern applications—particularly those involving high-dimensional data where these functional forms are often unknown, and both the drift and the diffusion must be learned directly from the observed trajectories. Statistical inference for SDEs has a rich history~\citep{KutoyantsBook2004}, with maximum-likelihood methods playing a central role when full trajectory data are available~\citep[Chapter $7$]{liptser2001statistics}. Recent advances have extended such methods to data-driven drift recovery \cite{guo2024learning}, but typically under restrictive noise assumptions, such as independence or constant variance.

In this work, we develop a noise-informed, trajectory-based learning framework for discovering the governing structures of SDEs directly from observational data. Unlike methods that estimate the drift alone or treat noise as a secondary effect, our approach embeds the noise process explicitly into the learning procedure and leverages information from the entire trajectory evolution, rather than focusing on isolated time points. This enables simultaneous recovery of both the drift $\vf$ and the noise covariance structure $\Sigma$, including scalar or matrix-valued forms and fully state-dependent, correlated noise. We conduct a systematic investigation of the method’s stability, accuracy, and computational efficiency across a variety of SDE and SPDE mdoels with different noise structures, demonstrating consistently superior performance in reconstructing complex stochastic dynamics.

The remainder of the paper is organized as follows.  We discuss the general SDE model which our learning is based on in Section~\ref{sec:model}.  Section~\ref{sec:Learning Framework} introduces the noise-informed likelihood formulation, the associated learning framework for recovering drift and noise, and a convergence theorem on learning drift. Section~\ref{sec:Examples} presents learning results on representative stochastic systems including high-dimensional interacting agent systems and stochastic heat equations, highlighting accuracy and robustness across diverse noise settings. Section~\ref{sec:supp_examples} discusses further supporting examples answering additional crucial aspects of our learning method . Section~\ref{sec:conclude} concludes with a discussion of the implications, limitations, and potential extensions of our approach.
\subsection{Related Works}\label{sec:Comparison}
System identification of the drift term from deterministic dynamics has been studied in many different scenarios, e.g. identification by enforcing sparsity such as SINDy~\citep{sindy}, neural network based methods such as NeuralODE~\citep{chen2018neural}, PINN~\citep{pinn2019} and autoencoder~\citep{xu2023modeling}, regression based~\cite{CS02}p, and high-dimensional reduction variational framework~\citep{LZTM2019}.  There are statistical methods which can be used to estimate the drift and noise terms using pointwise statistics.  SINDy for SDEs was also developed in~\citep{wanner2024higher}.  System identification using other dynamical properties are discussed in~\citep{HER2024, MIS2023}.

The observation data generated by SDEs can be treated as a time-series data with a mild assumption on the relationship between $\rvx_t$ and $\rvx_{t + \Delta t}$.  Various deep neural network architectures can be used to learn the drift term as well as predicting the trajectory data, using RNN, LSTM, and Transformers, see~\citep{liao2019learning, YANG2023279, wen2022transformers} for detailed discussion.

However, most of these methods use a regression type of loss function defined as follows
\[
\loss_{\gH}^{\text{Reg}}(\tilde\vf) = \E\Big[\frac{1}{T}\int_{t = 0}^T\norm{\tilde\vf(\rvx_t) - \frac{\dif\rvx_t}{\dt}}^2\, \dt\Big] + \lambda||\tilde\vf||_{Reg},
\]
where the constant $\lambda > 0$ is a penalty term assigned to control the weight of the regularization on the drift term with given a priori knowledge of the regularity of the drift.  Moreover, the derivative $\frac{\dif\rvx_t}{\dt}$ is loosely defined in the discrete sense (or weak sense). On the other hand, our likelihood induced loss of the form $\langle\tilde\vf, \Sigma^\dagger\tilde\vf\rangle\,\dt - 2\langle\tilde\vf, \Sigma^\dagger\dif\rvx_t\rangle$, is linked to the regression type loss through the expression 
\[
\norm{\tilde\vf - \frac{\dif\rvx_t}{\dt}}^2\, \dt = \norm{\tilde\vf}^2\, \dt - 2\langle\tilde\vf, \dif\rvx_t\rangle + \norm{\frac{\dif\rvx_t}{\dt}}^2\, \dt.
\]
The major difference comes in the re-scaling by the noise and our loss is a derivation from a negative-log likelihood, which guarantees the existence and uniqueness of minimizers.  See section \ref{sec:methods_comp} for a comparison of between these types of regression-based methods and our method.

Furthermore, special high-dim drift terms living on low-dim manifolds with constant noise is investigated in~\citep{lu2022}; such loss is similar to ours when $\sigma(\rvx) = \sigma > 0$.  In~\citep{guo2024learning}, a method constant correlated noise matrix is developed and studied.
\section{Model Equation}\label{sec:model}
Before introducing our learning framework for system identification from observed stochastic dynamics, we first establish the modeling setting and notation for the observational data.  Let $(\Omega,\mathcal{F},(\sF_{t})_{0\le t\le T},\mathbb{P})$ be a filtered probability space, for a fixed and  finite time horizon $T>0$.  As usual, the expectation operator with respect to $\mathbb{P}$ will be denoted by $\mathbb{E}_{\mathbb{P}}$ or simply $\mathbb{E}$. For random variables $X,Y$ we write $X\sim Y$, whenever $X,Y$ have the same distribution.  We consider governing equations for stochastic dynamics of the following form 
\begin{equation}\label{eq:original_model}
\dif \rvx_t = \vf(\rvx_t)\dif t + \sigma(\rvx_t) \dif \rvw_t, \quad \rvx_t, \rvw_t \in \R^D, 
\end{equation}
with some given initial condition $\rvx_0\sim \mu_0$, here $\vf: \R^D \rightarrow \R^D$ is the drift term, $\sigma: \R^D \rightarrow \R^{D \times D}$ is the diffusion coefficient.  Without Loss of Generality, we assume that $\sigma$ is symmetric positive definite (SPD), i.e., $\sigma^\intercal = \sigma$, $\rvx^\intercal\sigma\rvx \ge 0$ with $\rvx^\intercal\sigma\rvx = 0$ iff $\rvx = \mathbf{0}$.  Moreover, $\rvw$ represents a vector of independent standard Brownian motions. The covariance matrix of the SDE system is a symmetric positive definite matrix denoted by $\Sigma = \Sigma(\vx): \R^D \rightarrow \R^{D \times D}$ where $\Sigma = \sigma \sigma^\intercal$. We impose the following global regularity and growth conditions: there exist constants \(C_1, C_2 > 0\) such that for all \(\vx,\vy\in\R^{D}\)
\[
\begin{cases}
&\|\vf(\vx)-\vf(\vy)\|+\|\sigma(\vx)-\sigma(\vy)\|_{\mathrm{Fro}}
           \le C_1\|\vx-\vy\|, \\
&\|\vf(\vx)\|^{2}+\|\sigma(\vx)\|_{\mathrm{Fro}}^{2}
           \le C_2\bigl(1+\|\vx\|^{2}\bigr). 
\end{cases}
\]
Under these assumptions, \eqref{eq:original_model} admits a unique strong solution \(\{\rvx_t\}_{t\in[0,T]}\) adapted to the filtration \((\sF_{t})_{0\le t\le T}\) for every square-integrable initial condition \(\rvx_0\sim\mu_0\). 
\section{Learning Framework}\label{sec:Learning Framework}
We now introduce the methodology for learning the drift $\vf$ and the diffusion $\sigma$ terms of stochastic differential equations from observed trajectory data.  We assume continuous observation data $\{\rvx_t\}_{t \in [0, T]}$ for $\rvx_0 \sim \mu_0$, and that $\vf$ and $\sigma$ are the only unknowns.  We estimate these functions in two stages.
\subsection{Estimation of the Diffusion Term}
The diffusion coefficient $\sigma$ is first inferred using quadratic (co-)variation arguments.  For two scalar stochastic processes $\rx_t$ and $\ry_t$, the quadratic variation over time interval $[0, T]$ is defined by
\[
[\rx_t, \ry_t]_0^T = \lim_{|\Delta t_k| \to 0} \sum_{k= 1}^K (\rx(t_{k+1}) - \rx(t_k))(\ry(t_{k+1}) - \ry(t_k)),
\] 
where $\Delta t_k = t_{k + 1} - t_k$ and $\{0 = t_1 < t_2 < \cdots < t_K = T\}$ is a partition of the interval $[0, T]$.  For a vector stochastic process $\rvx_t = \begin{bmatrix}\rx_1(t), \rx_2(t), \ldots, \rx_D(t)\end{bmatrix}^\intercal $, the quadratic variation matrix $[\rvx,\rvx]_0^T$ has entries $[\rx_i(t), \rx_j(t)]_0^T$ for $i, j = 1, \ldots,D$.  Using such notation, the estimation of $\Sigma = \sigma\sigma^\intercal$ is the minimizer of the following loss function 
\begin{equation}\label{eq:sigma_loss}
\mathcal{E}_{\sigma}(\tilde \Sigma) = \E\Big[\big([\rvx_t, \rvx_t]_0^T - \int_{t=0}^T \tilde \Sigma(\rvx_t)\dif t\big)^2\Big].
\end{equation}
Since $\sigma$ is SPD, $\sigma = \sqrt{\Sigma}$ is uniquely defined. If $\Sigma$ is constant, then the estimation can be simplified to $\tilde \Sigma = \frac{1}{T}\mathbb{E}\big[[\rvx_t, \rvx_t]_0^T\big]$. Note that estimation of $\Sigma$ does not dependent on the drift function $\vf$.  
\subsection{Estimation of the Drift Term}
Once $\Sigma$ is obtained, we estimate $\vf$ by finding the minimizer to the following likelihood-based loss
\begin{equation}\label{eq:original_loss}
\loss_{\gH}(\tilde\vf) = \frac{1}{2}\E\Big[\int_{t = 0}^T\langle\tilde\vf(\rvx_t), \Sigma^{\dagger}(\rvx_t)\tilde\vf(\rvx_t)\rangle \, \dif t  - 2\langle\tilde\vf(\rvx_t), \Sigma^{\dagger}(\rvx_t)\dif \rvx_t\rangle\big)\Big],
\end{equation}
where $\tilde\vf \in \gH$ with $\gH$ being restricted to a convex and compact (w.r.t to $L_{\infty}$) function space determined by the observed data, $\langle\cdot,\cdot\rangle$ denotes the Euclidean inner product, and $\Sigma^\dagger$ is the pseudo-inverse of $\Sigma$, under our setting $\Sigma^\dagger = \Sigma^{-1}$.  The differential $\dif\rvx_t$ is approximated in practice by finite differences $\dif\rvx_t \approx \rvx_{t + \Delta t} - \rvx_t$.  This loss function arises from the Girsanov theorem and the Radon-Nikodym derivative for stochastic processes, see~\citep[Chpater~7]{liptser2001statistics} and Section~\ref{sec:theory} for details.  
\subsection{Derivation of the Loss for the Drift}\label{sec:theory}
We discuss the theoretical foundation of our methods in this section. Consider two It\^o processes defined over measurable space $(\Omega, \mathcal{F})$ and let $\mathbb{P}_X$, $\mathbb{P}_Y$ be probability measures corresponding to processes $\rvx$ and $\rvy$, where
\[
\begin{aligned}
  \dif \rvx_t &= \vf(\rvx_t)\dif t + \sigma(\rvx_t)\dif \rvw_t, \\
  \dif \rvy_t &= \vg(\rvy_t)\dif t + \sigma(\rvy_t)\dif \rvw_t, \quad \rvy_0 = \rvx_0,
\end{aligned}
\]
satisfying all assumptions in~\citep[Theorem~7.18]{liptser2001statistics} and its following corollary. Then, the Radon-Nikodym derivative, or the likelihood ratio, takes the form
\begin{equation}
\frac{\dif \mathbb{P}_X}{\dif \mathbb{P}_Y}(\rvy) = \exp \Big( \int_0^T \inp{(\vf_t - \vg_t, \Sigma^{\dagger}\dif \rvy_t} -\frac{1}{2} \int_0^T \inp{(\vf_t - \vg_t), \Sigma_t^{\dagger} (\vf_t + \vg_t}\dif t\Big),
\end{equation}
where $\vf_t = \vf(\ry_t)$, $\vg_t = \vg(\rvy_t)$, and $\Sigma_t = \Sigma(\rvy_t)$. Here let us assume that the observations are $\{\rvx_t\}_{t \in [0, T]}$. In view of the assumption of~\citep[Theorem~7.18]{liptser2001statistics}, the  $n$-dimensional adapted process $\Theta = \sigma^{\dagger}(\vf(\rvx_t) - \vg(\rvx_t))$  is such that $\int_0^T \norm{\Theta}^2\dif t < \infty$. By Girsanov theorem, $\widetilde{\rvw_t} = \rvw_t + \int_0^T \Theta_s \dif s$ is an $n$-dimensional standard Brownian motion under probability measure $\mathbb{P}_Y$. Hence,  $\dif\rvx_t = \vf(\rvx_t)\dif t + \sigma(\rvx_t)(\dif \tilde{\rvw_t} - \Theta_t \dif t) = \vg(\rvx_t)\dif t + \sigma(\rvx_t)\dif \tilde{\rvw_t}$. For convenience,  we take $\vg=0$, in which case $\rvx_t$ becomes a Brownian process under $\mathbb{P}_Y$. Therefore $\mathbb{P}_Y(\{\rvx_t\}_{t \in [0, T]} | \vf)$ is now independent from $\vf$ since $\rvx_t$ has no drift term under $\mathbb{P}_Y$. Putting such likelihood under the negative-log function, we arrive at our first loss as
\[
\loss_T(\tilde\vf) = -\ln{L(\vf|\{\rvx_t\}_{t \in [0, T]})} = \int_0^T\big(\vf(\rvx_t)^\intercal \Sigma^{\dagger} f(\rvx_t)\dif t - 2\vf(\rvx_t)^\intercal \Sigma^{\dagger}\dif\rvx_t\big).
\]
Here such loss function is used to handle observation data from one long trajectory (i.e. observed over large time), and it will be effective especially for ergodic systems.  Moreover, we also consider the situation where multiple medium (or short-burst) trajectories with different initial conditions are observed, then we derive our loss function as the expectation (over trajectories with different initial conditions) of the negative-log-likelihood function as
\[
\loss(\tilde\vf) = \E\big[-\ln{L(\vf|\{\rvx_t\}_{t \in [0, T]})}\big] = \frac{1}{2}\E\big[\int_0^T\big(\vf(\rvx_t)^\intercal \Sigma^{\dagger} f(\rvx_t)\dif t - 2\vf(\rvx_t)^\intercal \Sigma^{\dagger}\dif\rvx_t\big)\big].
\]
\subsection{Convergence Theorem}
We present the following convergence results in a theorem.
\begin{theorem}
Given the continuous-time i.i.d trajectory data $\{\rvx^m_t\}_{m = 1}^M$ for $t \in [0, T]$ and each $\rvx_t^m$ generated by \eqref{eq:original_model}, we define an estimator to $\vf$ through minimizing the following loss
\[
\loss_M(\tilde\vf) = \frac{1}{2M}\sum_{m = 1}^M\big(\int_0^T\inp{\tilde\vf_t^m, \,(\Sigma^m_t)^{-1}\tilde\vf_t^m} \, \dt - 2\int_0^T\inp{\tilde\vf_t^m, \,(\Sigma^m_t)^{-1}\dif\rvx_t^m}\big),
\]
where $\tilde\vf_t^m = \tilde\vf(\rvx_t^m)$, $\Sigma_t^m = \Sigma(\rvx_t^m)$, and $\tilde\vf \in \sH$ with $\sH$ being convex and compact (w.r.t to $L^2$-norm).  When $\sH$ is finite dimensional, i.e., $n = \dim(\sH) < \infty$, and $\vf \in \sH$, then the estimator, given as $\hat\vf_{M} = \argmin_{\tilde\vf \in \sH}\loss_M(\tilde\vf)$, has the following properties: $\hat\vf_M \xrightarrow{P} \vf$ (consistency) and $\sqrt{M}(\hat\vf_M - \vf) \xrightarrow{D} \normal(\vzero, \mB^{-1})$ (Asymptotic normality). Here $\mB = \E[\int_0^T\Psi_t^\intercal\Sigma_t^{-1}\Psi_t \, dt]$, where 
\[
\Psi_t = \Psi(\rvx_t) = \begin{bmatrix} \psivec_1(\rvx_t) & \cdots & \psivec_n(\rvx_t) \end{bmatrix} \in \R^{D \times n}.
\]
with $\{\psivec_1, \psivec_2, \cdots, \psivec_n\}$ being a basis of $\sH$ where each $\psivec_{\eta}: \R^D \rightarrow \R^D$.  Notice that $\mB$ is SPD. 
\end{theorem}
A detailed proof of this theorem is presented in Appendix \ref{append:proof}.
\subsection{Deep Learning for High-dim Functions}\label{sec:DL_learn}
In learning high dimensional $\vf$ and $\sigma$, we can employ the deep learning architecture, with one neural network for learning $\vf$ and the other for $\sigma$.  The learning of $\vf$ is rather straightforward, since the loss is well-defined for deep learning and simply changing the functional space to be a space of neural networks.  We will discuss the learning of $\sigma$ in details. Let $G:\R^{D}\to\R^{D(D+1)/2}$ be a neural network with outputs arranged as $\{\vu_{ij}(\vx)\}_{1\le j\le i\le D}$.   Since $\Sigma$ is SPD, the Cholesky decomposition on $\Sigma$ gives $\Sigma(\vx):=L(\vx)\,L(\vx)^\intercal$ where $L$ is a lower-triangular matrix with positive diagonal entries.  Therefore we can learn a lower–triangular mapping $\tilde{L}:\R^{D}\to\R^{D\times D}$ by
\[
\big(\tilde{L}(\vx)\big)_{ij} = \begin{cases}
h\big(\vu_{ii}(\vx)\big) & \text{if } i = j \\
\vu_{ij}(\vx) & \text{if } i > j \\
0 & \text{if } i < j
\end{cases},
\]
where $h:\R\to(0,\infty)$ is some chosen function to enforce positivity on the main diagonal. Hence, we define the model $\tilde \Sigma(\vx):=\tilde{L}(\vx)\,\tilde{L}(\vx)^\intercal \approx \Sigma(\vx)$. Given $M$ trajectories, set $Y^{m}_l\coloneqq\frac{\Delta\rvx^{m}_l\big(\Delta\rvx^{m}_l\big)^\intercal}{\Delta t}$.  We learn the estimator by minimizing the Frobenius mean squared difference between $Y_l$ and $\tilde \Sigma(\rvx_{l})$ over all observed trajectories:
\begin{equation}\label{eq:QVMM}
\mathcal \loss(\tilde \Sigma)
= \frac{1}{M}\sum_{m=1}^{M} \sum_{l=0}^{L-1}
\Big\|\,Y^{m}_l-\tilde \Sigma\big(\rvx^{m}_{l}\big)\,\Big\|_{\mathrm F}^{2}.
\end{equation}

If $\sigma$ is a full matrix, we use the matrix-square-root function to obtain $\sigma = \sqrt{\Sigma}$.  If $\Sigma(\rvx)$ is diagonal for all $\vx$, i.e., $\Sigma(\vx)=\mathrm{diag}\big(\Sigma_{11}(\vx),\dots,\Sigma_{DD}(\vx)\big)$, then we will learn each diagonal entry by a single-output positive network. Writing $Y^{m}_{l,ii}=\frac{\big(\Delta \rvx^{(m,i)}_l\big)^2}{\Delta t}$, where \(\rvx^{(m,i)}_l \) represents the $i^{th}$ entry of \(\rvx_l^m\) and the loss function can be decoupled and become $\mathcal \loss(\tilde \Sigma_{ii})
=\frac{1}{M}\sum_{m=1}^{M}\sum_{l=0}^{L-1}\Big(\,Y^{m}_{l,ii}-\tilde\Sigma_{ii}\big(\rvx^{m}_{l}\big)\,\Big)^{2}$.Hence $\hat\sigma_{ii}(\vx)=\sqrt{\,\hat\Sigma_{ii}(\vx)\,}$.
\subsection{Performance Measures}\label{subsec:Performance Measures}
In order to properly gauge the accuracy of our learning estimators, we provide three different performance measures of our estimated drift.  First, if we have access to original drift function $\vf$, then we will use the following error to compute the difference between $\hat\vf$ (our estimator) to $\vf$ with the following norm
\begin{equation}\label{eq:rho_norm}
    \norm{\vf - \hat\vf}_{L^2(\rho)}^2 = \int_{\R^d}\norm{\vf(\rvx) - \hat\vf(\rvx)}_{\ell^2(\R^D)}^2 \, \dif \rho(\rvx),
\end{equation}
where the weighted measure $\rho$,  defined on $\R^D$, is $\rho(\rvx) = \E\Big[\frac{1}{T}\int_{t = 0}^T\delta_{\rvx_t}(\rvx)\Big]$.  Here $\rvx_t$ evolves from $\rvx_0$ by \eqref{eq:original_model}. The norm given by \ref{eq:rho_norm} is useful only from the theoretical perspective, e.g. showing convergence.  Under normal circumstances, $\vf$ is most likely non-accessible.  Thus we look at a performance measure that compares the difference between $\rmX(\vf, \rvx_0, T) = \{\rvx_t\}_{t \in [0, T]}$ (the observed trajectory that evolves from $\rvx_0 \sim \mu_0$ with the unknown $\vf$) and $\hat\rmX(\hat\vf, \rvx_0, T) = \{\hat\rvx_t\}_{t \in [0, T]}$ (the estimated trajectory that evolves from the same $\rvx_0$ with the learned $\hat\vf$ and driven by the same realized random noise as used by the original dynamics).  Then, the difference between the two trajectories is measured as follows
\begin{equation}\label{eq:traj_norm}
\norm{\rmX - \hat\rmX} = \E\Big[\frac{1}{T}\int_{t = 0}^T\norm{\rvx_t - \hat\rvx_t}_{\ell^2(\R^D)}^2 \, \dif t\Big].
\end{equation}
However, comparing two sets of trajectories (even with the same initial condition) on the same random noise is not realistic. Therefore we compare the distribution of the trajectories over different initial conditions and different noise at the same time instances using the Wasserstein distance at any given time $t \in [0, T]$. Let $\mu^M_t$ be the empirical distribution at time $t$ for the simulation under $\vf$ with $M$ trajectories, and $\hat{\mu}^M_t$ be the empirical distribution at time $t$ for the simulation with $M$ trajectories under $\hat{\vf}$, where $\mu^M_t = \frac{1}{M} \sum_{i=1}^M \delta_{\rvx^{(i)}(t)}$, $\hat{\mu}^M_t = \frac{1}{M} \sum_{i=1}^M \delta_{\hat{\rvx}^{(i)}(t)}$. Then the Wasserstein distance of order two between $\mu^M_t$ and $\hat{\mu}^M_t$ is defined  as
\begin{equation}\label{eq:wass_dis}
\mathcal{W}_2(\mu^M_t, \hat{\mu}^M_t \,|\, \mu_0)= \left( \inf_{\pi \in \Pi(\mu^M_t, \hat{\mu}^M_t \,|\, \mu_0)} \int_{\R^D \times \R^D} \|x - y\|^2 \, \dif \pi(x, y) \right)^{1/2}.
\end{equation}
Here, $\Pi(\mu^M_t, \hat{\mu}^M_t \,|\, \mu_0)$ is the set of all joint distributions on $\R^D \times \R^D$ with marginals $\mu^M_t$ and $\hat{\mu}^M_t$, and with the additional constraint that the joint distribution must be consistent with the initial distribution of $\rvx_0$ following $\mu_0$.
\section{Examples}\label{sec:Examples}
We demonstrate the application of our trajectory-based method for estimating drift functions and noise structures through a set of representative stochastic systems.  The examples are organized into three categories: (i) a benchmark low-dimensional SDE with state-dependent diffusion, (ii) high-dimensional interacting particle systems (IPS) with structured drift, and (iii) stochastic partial differential equations (SPDEs) after a Galerkin/spectral discretization, where the effective dimension can become large.

\subsection{Example: Benchmark Model}\label{sec:example_heston}
\begin{figure}[h!]
    \centering
    \includegraphics[width=0.9\linewidth]{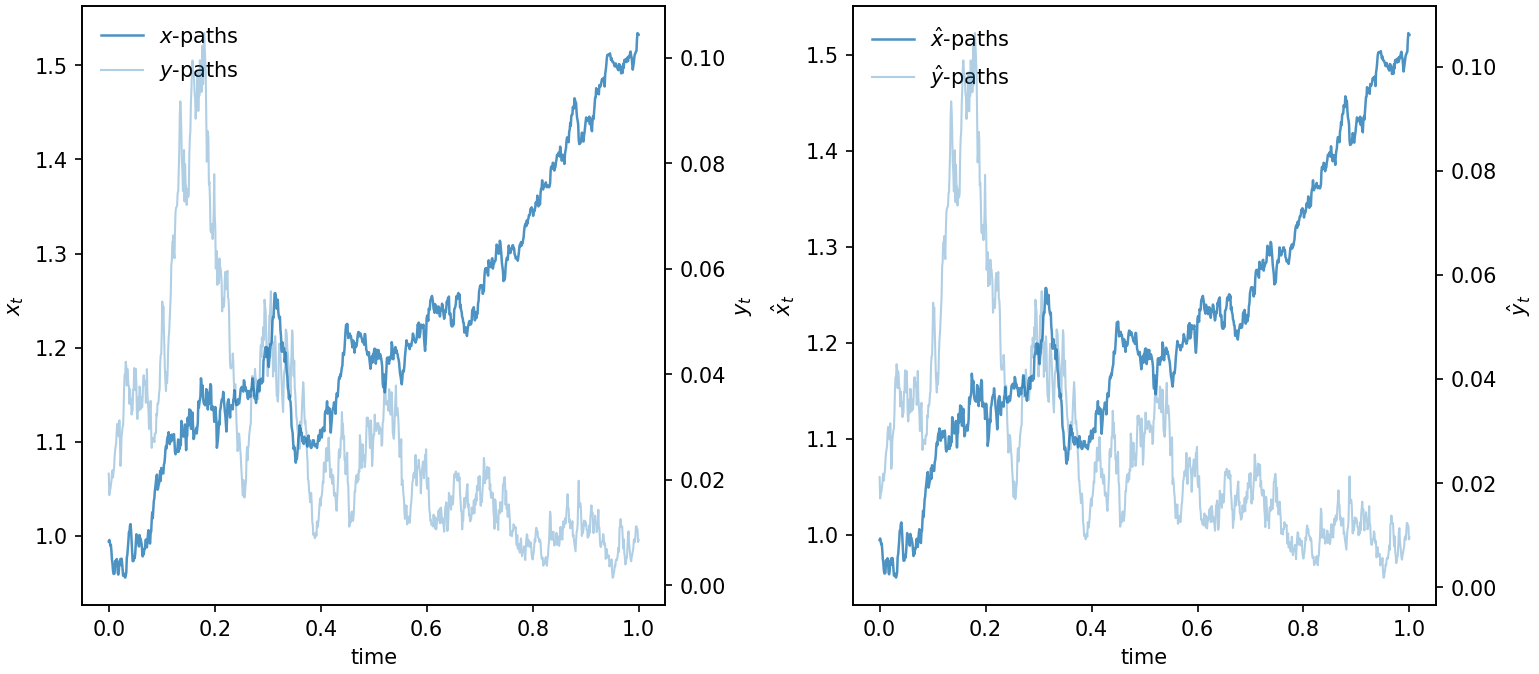}
    \caption{Trajectory comparison with matched noise \( \dif\rb_t \). Left: true simulated paths $(\rx_t,\ry_t)$ under the benchmark parameters. Right: re–simulated paths using the learned drift $\hat \rf$ and diffusion $\hat\sigma$, driven by the same \( (\dif\rb_t^x, \dif\rb_t^y ) \).}
    \label{fig:benchmark_3}
\end{figure}

We consider an SDE model with state dependent noise matrix, as follows
\[
\begin{cases}
    \dif\rx_t &= C_1\rx_t\dt +\sqrt{\ry_t}\rx_t\dif\rb_t^x \\
    \dif\ry_t &= C_2(C_3 - \ry_t)\dt + C_4\sqrt{\ry_t}\dif\rb_t^y,
\end{cases}
\]
where $(\rx_t, \ry_t)$ is the pair of state-variables, $(\rb_t^x, \rb_t^y)$ are standard Brownian motion, the constants $C_1, C_2 , C_3, C_4 > 0$ are model parameters.  If $2C_2C_3 > C_4^2$, then $\ry_t$ remains strictly positive.  We use this benchmarking model to test the effectiveness of our learning framework on identifying the SDE without any knowledge of the noise and drift terms.
\begin{figure}[h!]
  \begin{subfigure}[b]{0.45\textwidth}
  \centering
  \includegraphics[width=\textwidth]{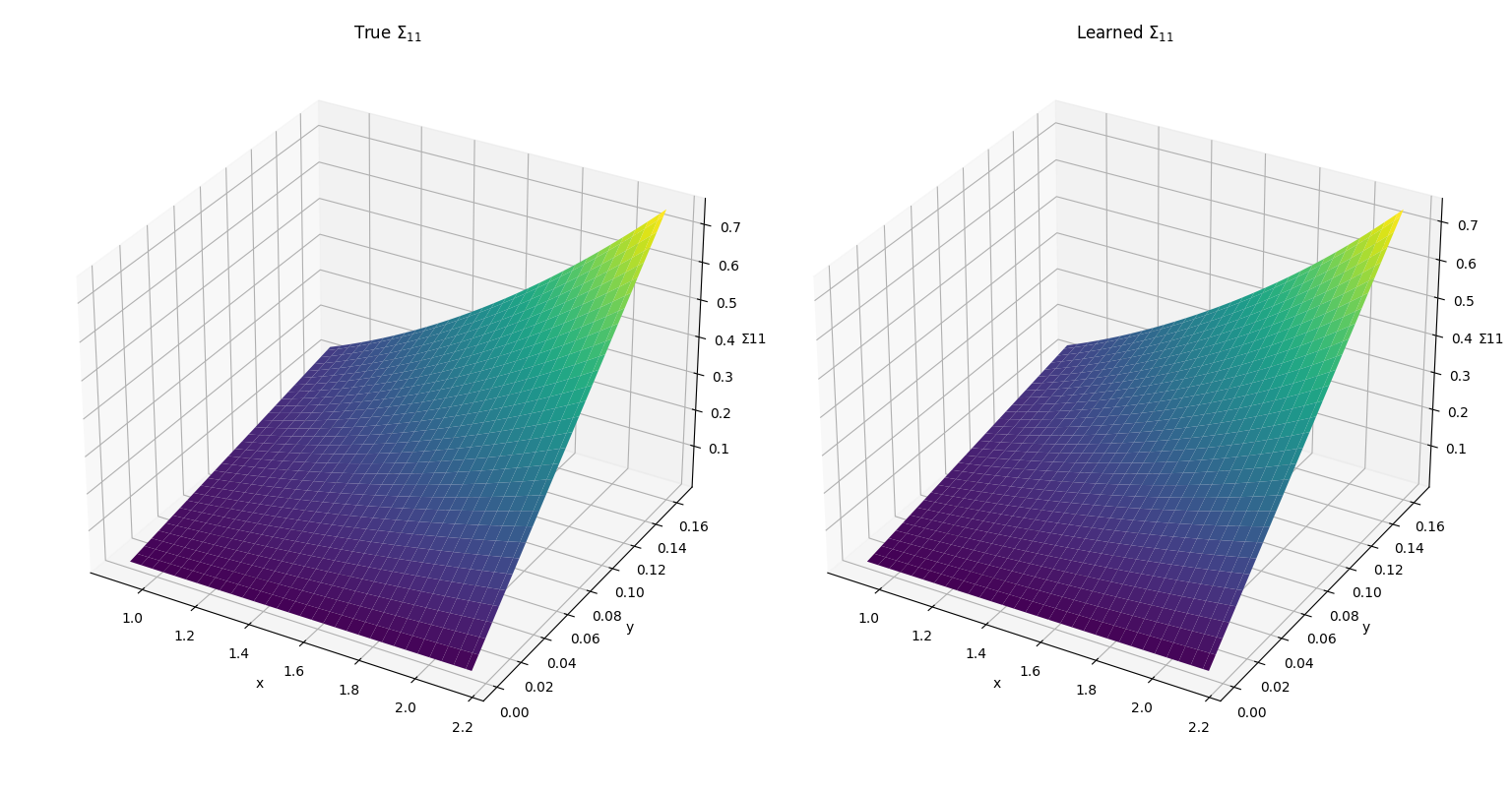}
  \caption{Comparison of true $\Sigma_{11}$ vs learned $\hat \Sigma_{11}$. }
  \label{fig:benchmark_sigma_1}
  \end{subfigure}
  \hfill
  \begin{subfigure}[b]{0.45\textwidth}
  \centering
  \includegraphics[width=\textwidth]{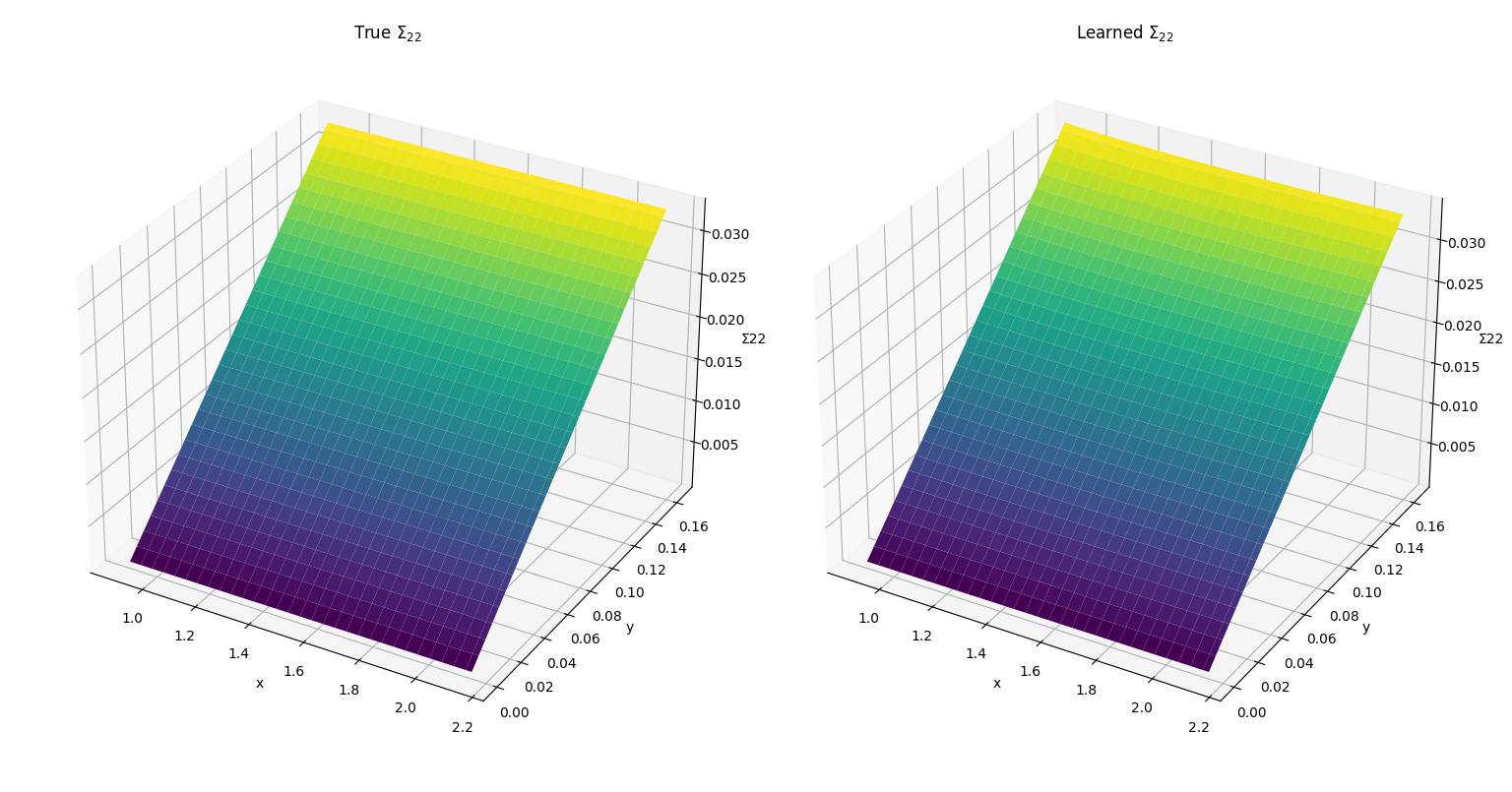}
  \caption{Comparison of true $\Sigma_{22}$ vs learned $\hat \Sigma_{22}$.}
  \label{fig:benchmark_sigma_2}
  \end{subfigure}
\caption{Benchmark model: $\Sigma$ vs $\hat\Sigma$.}
\label{fig:benchmark_1}
\end{figure}
We evaluated our learning method on the benchmark model. Trajectories were simulated using the parameters $C_1=0.5$, $C_2=3.0$, $C_3=0.04$, and $C_4=0.45$. Both the drift function $f(x, y) = [f_1(x), f_2(y)]$, where $f_1(x) = C_1x$ and $f_2(y) = C_2(C_3-y)$, and the diffusion matrix $\sigma(x, y) = \begin{bmatrix} \sqrt{y}x & 0 \\ 0 & C_4\sqrt{y} \end{bmatrix}$ were learned using the neural network method described in the previous section \ref{sec:DL_learn}. 
\begin{figure}[H]
\begin{subfigure}[b]{0.45\textwidth}
  \centering
  \includegraphics[width=\textwidth]{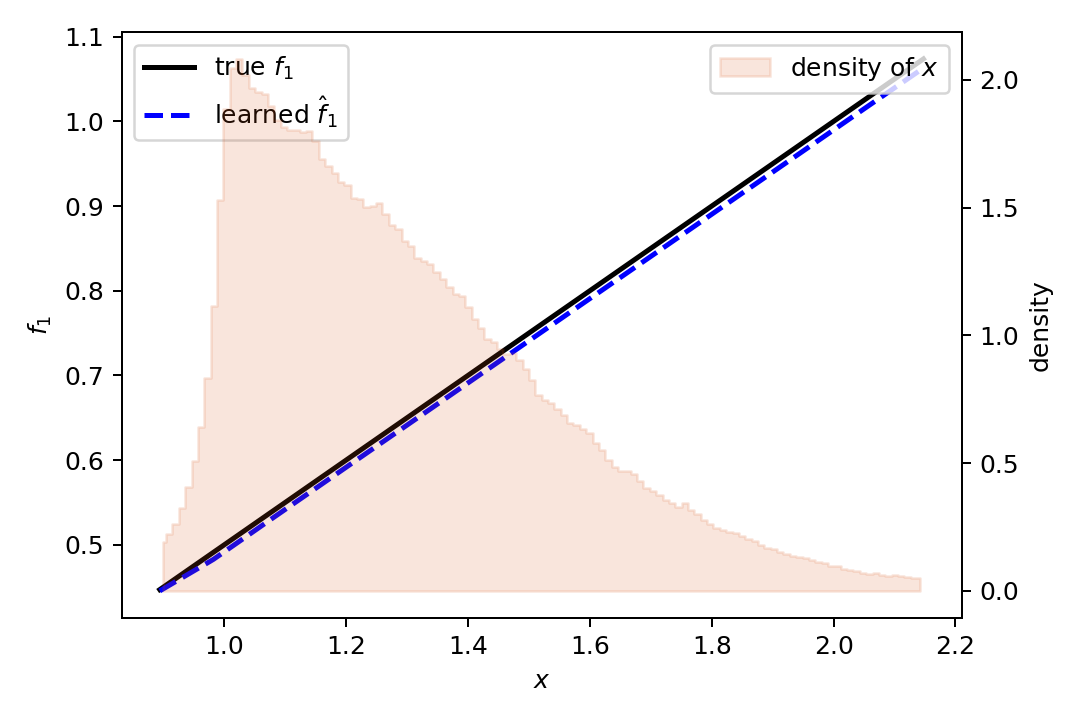}
  \caption{True $f_1$ vs learned $\hat{f}_1$.}
  \label{fig:benchmark_drift_1}
  \end{subfigure}
  \hfill
  \begin{subfigure}[b]{0.45\textwidth}
  \centering
  \includegraphics[width=\textwidth]{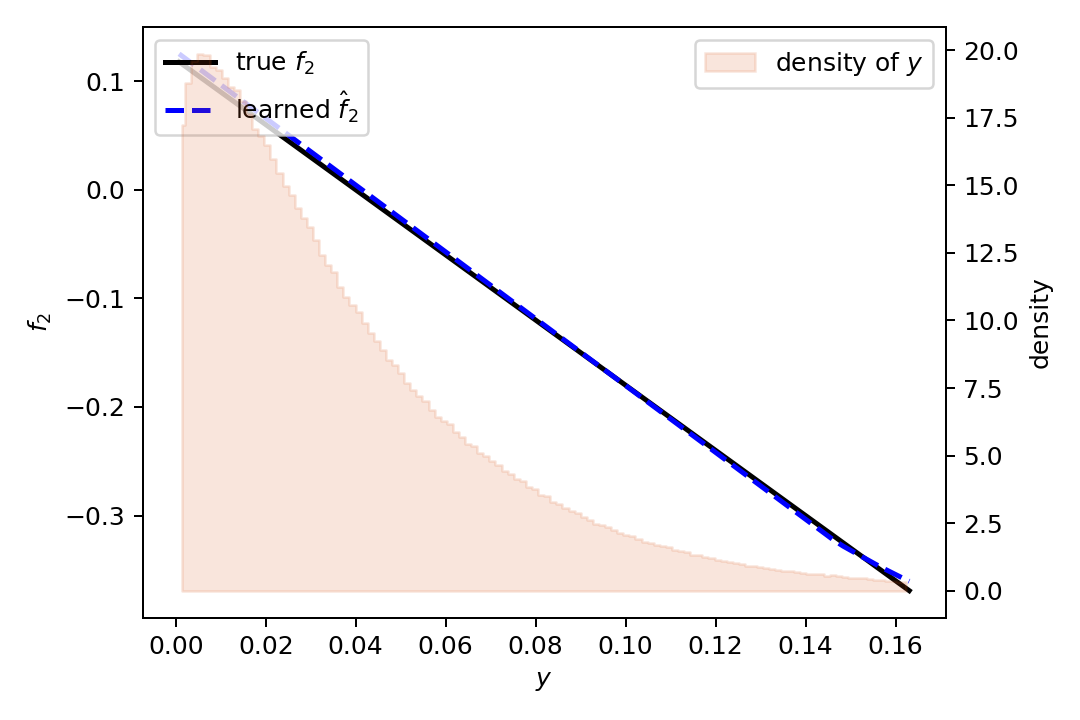}
  \caption{True $f_2$ vs learned $\hat{f}_2$.}
  \label{fig:benchmark_drift_2}
  \end{subfigure}
\caption{$\vf$ vs $\hat\vf$ with empirical distribution of $\rvx_t$ is shown in the background.}
\label{fig:benchmark_2}
\end{figure}
The results are shown in Figure \ref{fig:benchmark_1}, \ref{fig:benchmark_2} and \ref{fig:benchmark_3}.
By using deep neural networks as the underlying function spaces, one can easily infer those multi-variate drift and noise functions, without specifying the actual form of the functions.   

\subsection{Example: Interacting Particle Systems (IPS)}\label{sec:example_agents}
\begin{figure}[ht!]
  \begin{subfigure}[b]{0.45\textwidth}
  \centering
  \includegraphics[width=\textwidth]{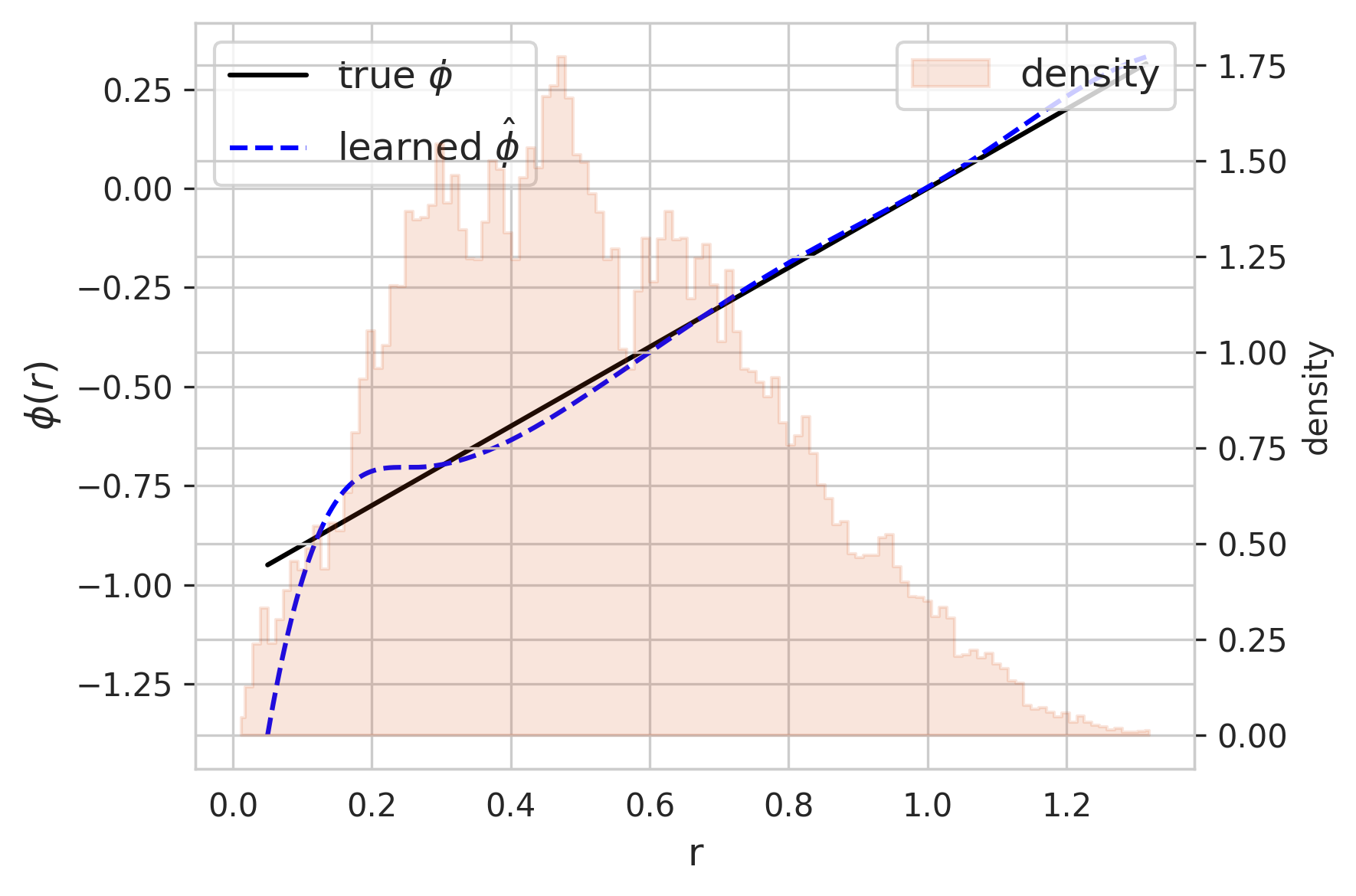}
  \caption{Case (I): Err = $0.02$.}
  \label{fig:60D_phi_1}
  \end{subfigure}
  \hfill
  \begin{subfigure}[b]{0.45\textwidth}
  \centering
  \includegraphics[width=\textwidth]{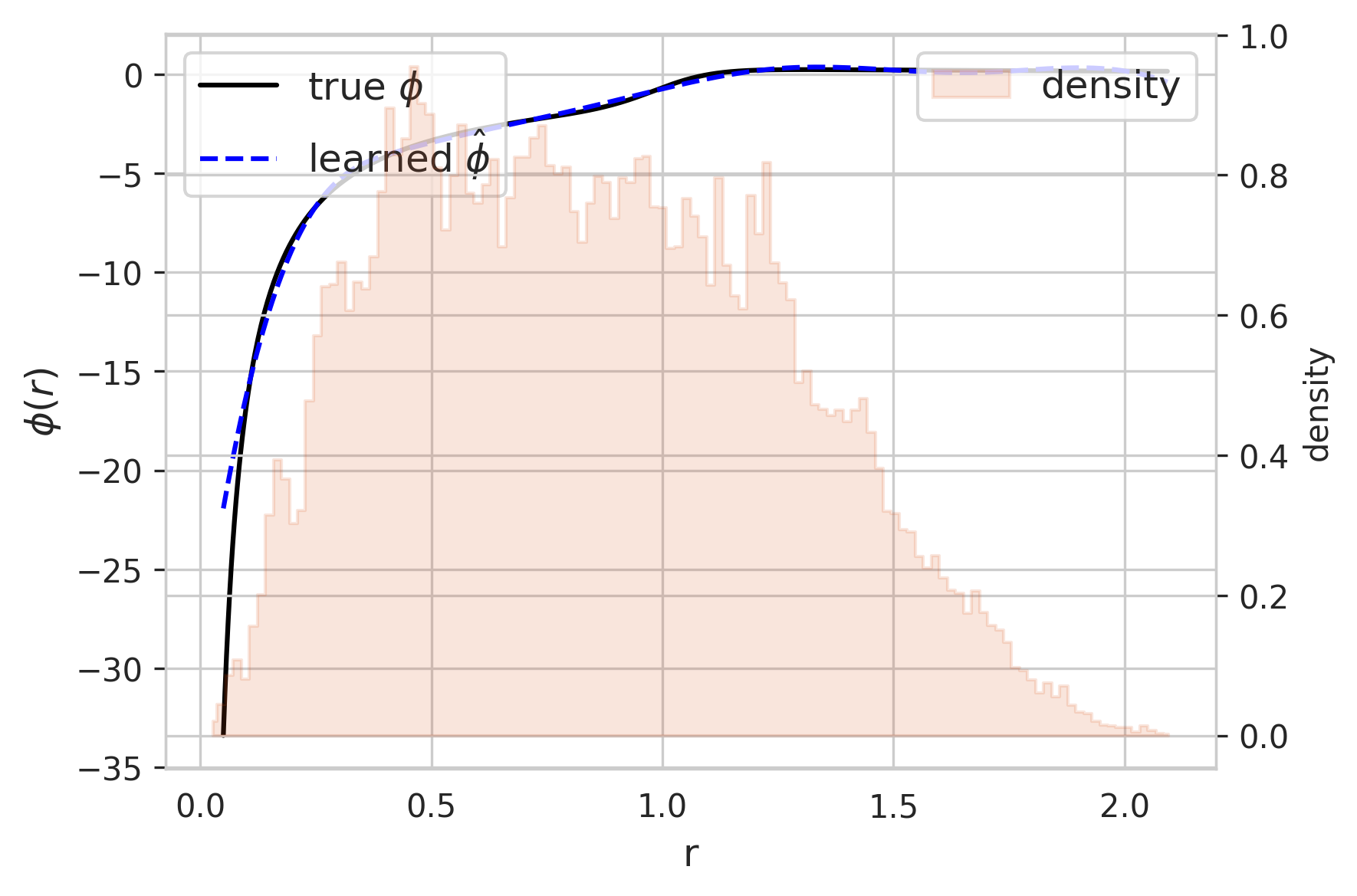}
  \caption{Case (II): Err = $0.14$.}
  \label{fig:60D_phi_2}
  \end{subfigure}
\caption{True $\phi$ vs learned $\hat\phi$; Empirical density of $r$ shown in the background.}
\label{fig:60D case phi}
\end{figure}
We consider a high dimensional SDE case where the drift term has a special structure. Such special structure will allow us to learn the high-dimensional SDE more effectively through an innate dimension reduction approach.  This high dimensional SDE case is a presentation of an interacting particle system.  Learning of such systems without stochastic noise terms had been investigated in~\citep{LZTM2019, ZMM2020, MMQZ2021, FMMZ2022, FM2023}.  We consider such system with correlated and state-dependent stochastic noise, i.e. for a system of $N$ particles, where each particle is associated with a state vector $\rvx_i \in \R^{d}$.  The particles' states are governed by the following system of SDEs
\[
\dif \rvx_i(t) = \frac{1}{N}\sum_{j = 1, j \neq i}^N\phi(\norm{\rvx_j(t) - \rvx_i(t)})(\rvx_j(t) - \rvx_i(t))\dif t + \sigma^{\rx}(\rvx_i(t))\dif \rvw(t), \quad i = 1, \ldots, N.
\]
Here $\phi:\R^+ \rightarrow \R$ is an interaction kernel that governs how particle $j$ influences the behavior of particle $i$, and $\sigma^{\rx}: \R^{d} \rightarrow \R^{d \times d}$ is a symmetric positive definite matrix that represents the noise strength and correlation.  
\begin{figure}[ht!]
  \begin{subfigure}[b]{0.45\textwidth}
  \centering
  \includegraphics[width=\textwidth]{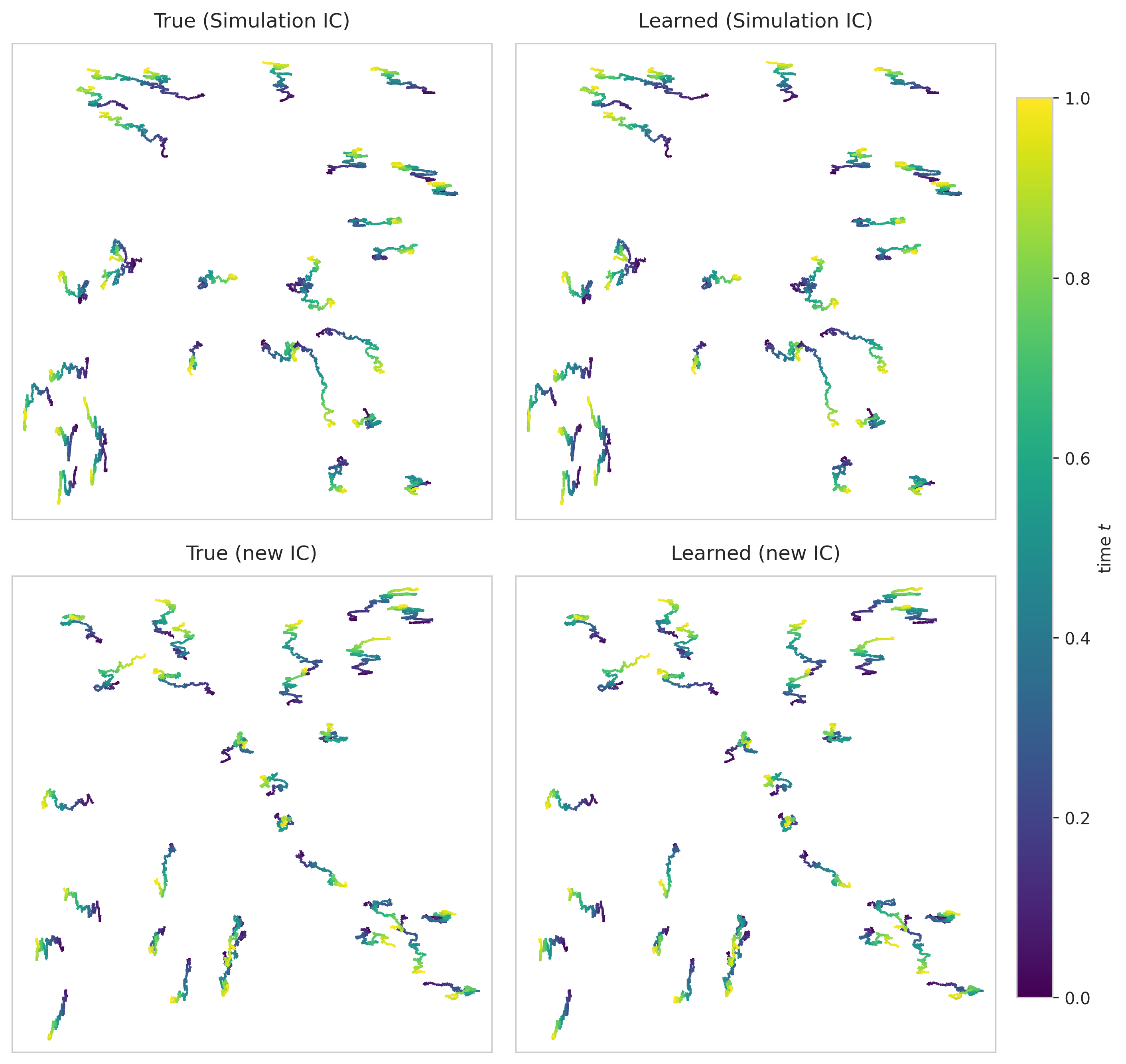}
  \caption{Case (I).}
  \label{fig:60D_trajs_1}
  \end{subfigure}
  \hfill
  \begin{subfigure}[b]{0.45\textwidth}
  \centering
  \includegraphics[width=\textwidth]{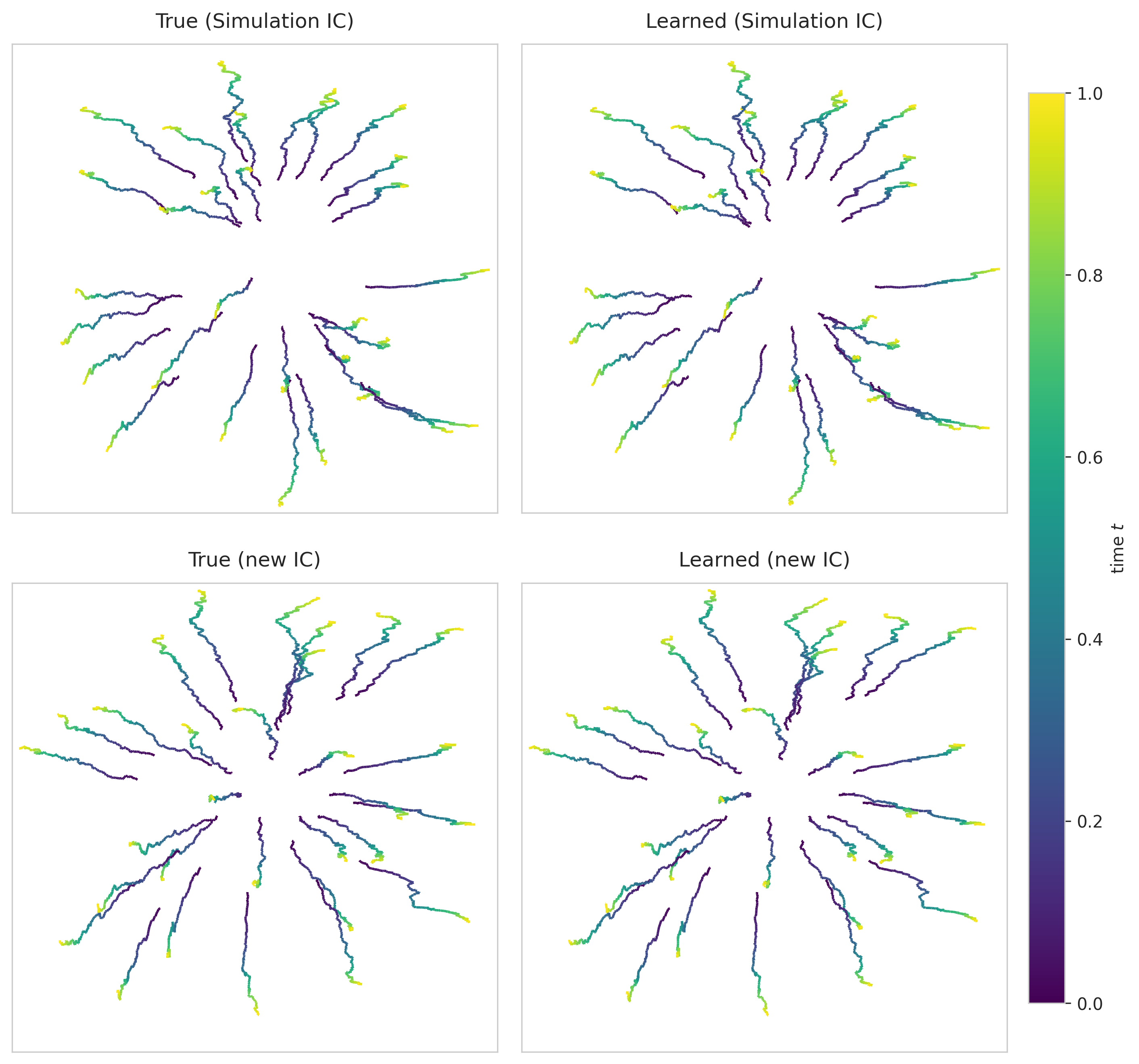}
  \caption{Case (II).}
  \label{fig:60D_trajs_2}
  \end{subfigure}
\caption{True $\rvx$ vs learned $\hat\rvx$ under the same noise. Top row: evolution from the same training IC. Bottom row: evolution from a new IC.}
\label{fig:60D case trajs}
\end{figure}
We test two interaction kernels
\[
\begin{aligned}
\text{Case (I)}: &\quad \phi(r)=r-1, \\
\text{Case (II)}: &\quad \phi(r)=-\frac{\tanh\!\big(8(1-r)\big)+0.67}{r}.
\end{aligned}
\]
The diffusion is shared across particles, diagonal, and state–dependent, i.e., $\sigma^{\rx}(\rvx_i(t))=\operatorname{diag}\big(\sigma^{\rx}_{11}(\rvx_i(t)),\,\sigma^{\rx}_{22}(\rvx_i(t))\big)$ with
\[
\begin{cases}
\sigma^{\rx}_{11}(\rvx_i(t))&=0.08\,\sin^{2}\!\big(\|\rvx_i(t)\|\big)+\varepsilon,\\
\sigma^{\rx}_{22}(\rvx_i(t))&=0.06\,\cos^{2}\!\big(\|\rvx_i(t)\|\big)+\varepsilon,
\end{cases} \quad \varepsilon=0.01.
\]
We run two experiments to justify our method. We take \(N=30\) particles in \(\R^{d}\) with \(d=2\) (so \(D=Nd=60\)), time horizon \(T=1\), step size \(\Delta t=0.001\), and \(M=100\) i.i.d.\ trajectories. The initial distributions are i.i.d.\ \(\rvx_0\sim \mathrm{Unif}([0,1]^{d})\) for each particle. Simulation uses Euler–Maruyama method. In estimating $\sigma$, following the general implementation of \(\sigma\) mentioned in ~\ref{subsec: Implementation}, for the diagonal case we learn each diagonal entry independently.
\begin{figure}[ht!]
  \begin{subfigure}[b]{0.45\textwidth}
  \centering
  \includegraphics[width=\textwidth]{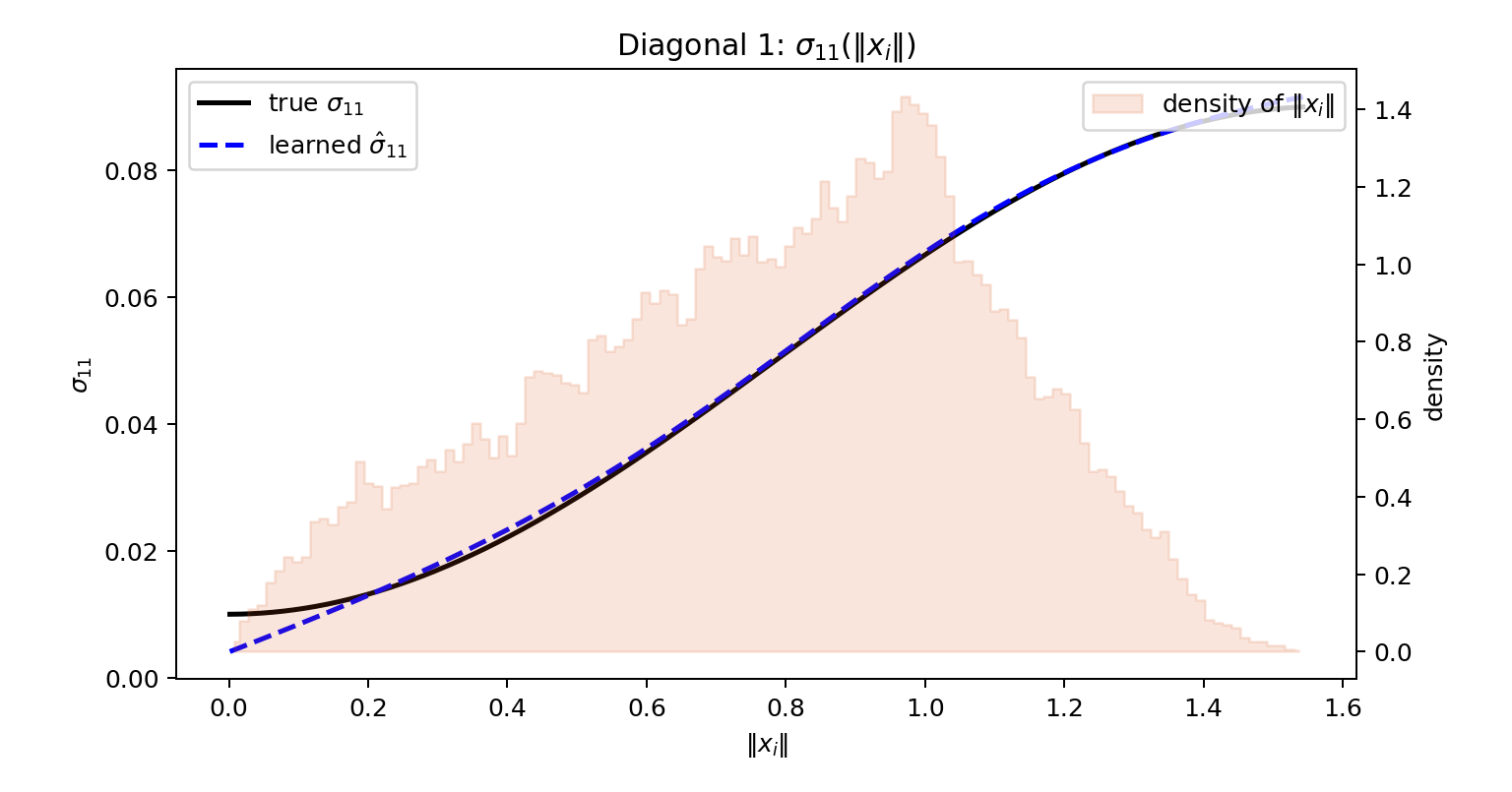}
  \caption{Case (I): Err = $0.013$.}
  \label{fig:Agent_sigma_1_11}
  \end{subfigure}
\hfill
  \begin{subfigure}[b]{0.45\textwidth}
  \centering
  \includegraphics[width=\textwidth]{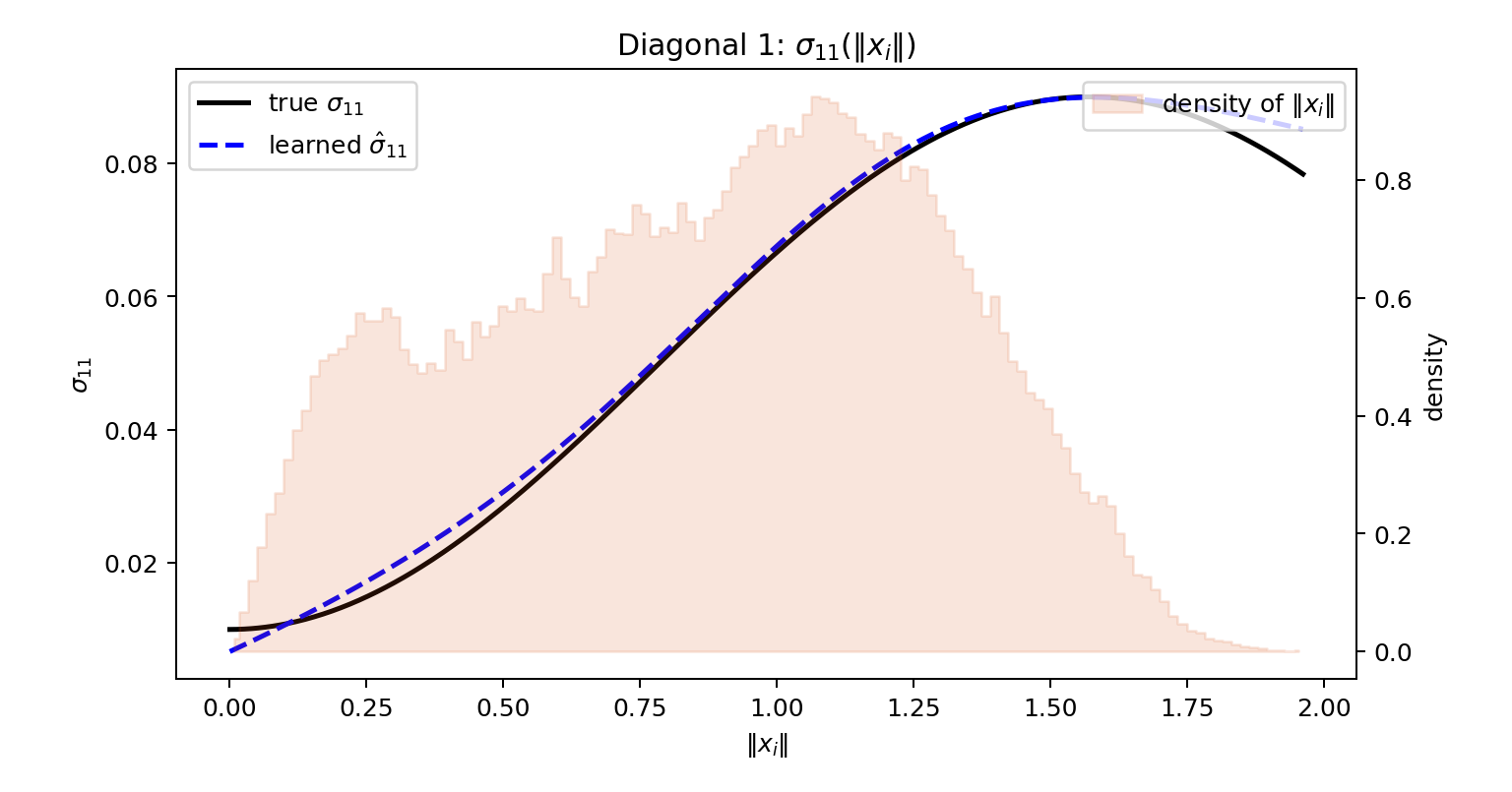}
  \caption{Case (II): Err = $0.023$}
  \label{fig:Agent_sigma_1_22}
  \end{subfigure}
\\
  \begin{subfigure}[b]{0.45\textwidth}
  \centering
  \includegraphics[width=\textwidth]{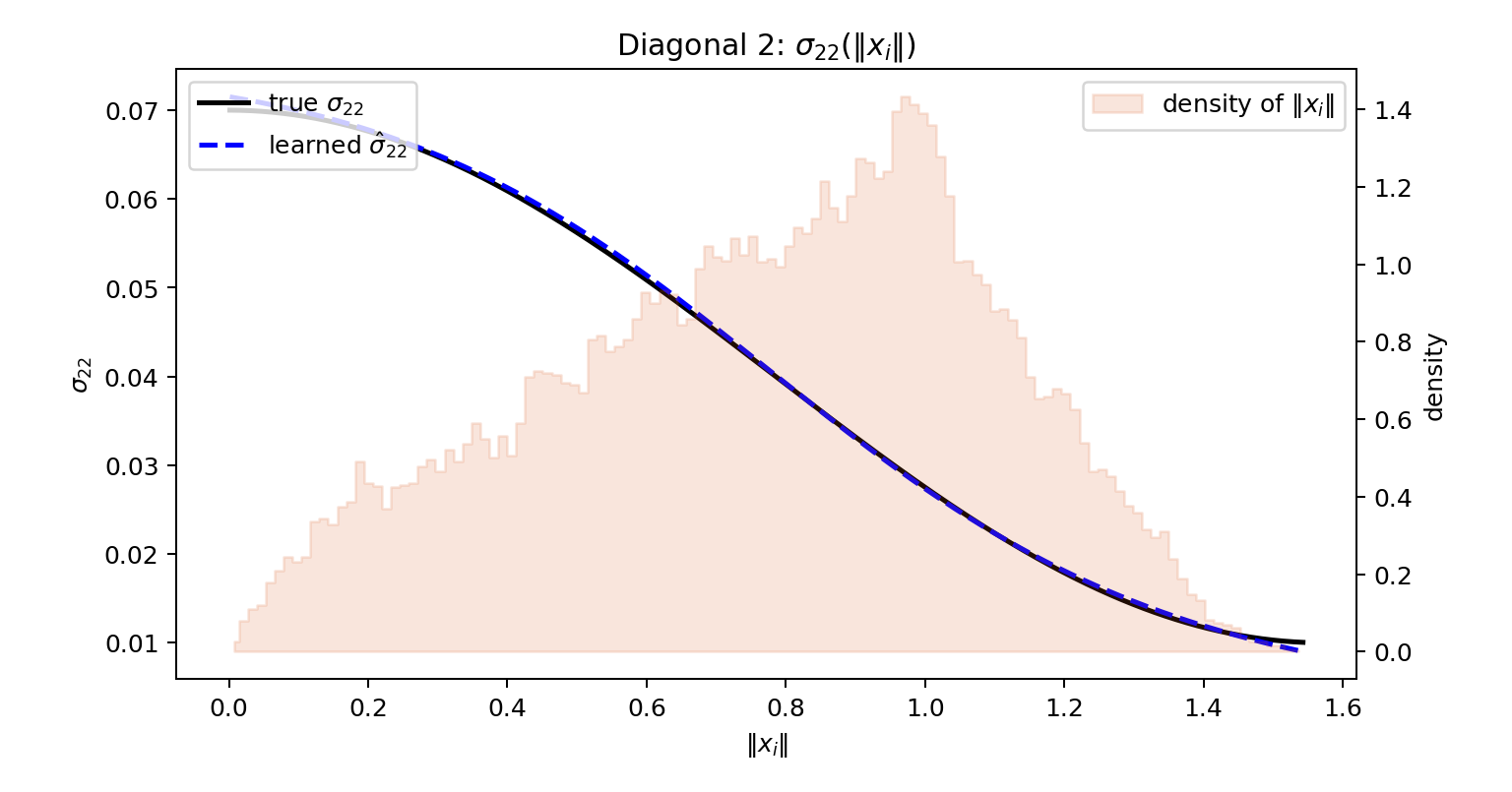}
  \caption{Case (I): Err = $0.007$.}
  \label{fig:Agent_sigma_2_11}
  \end{subfigure}
\hfill
  \begin{subfigure}[b]{0.45\textwidth}
  \centering
  \includegraphics[width=\textwidth]{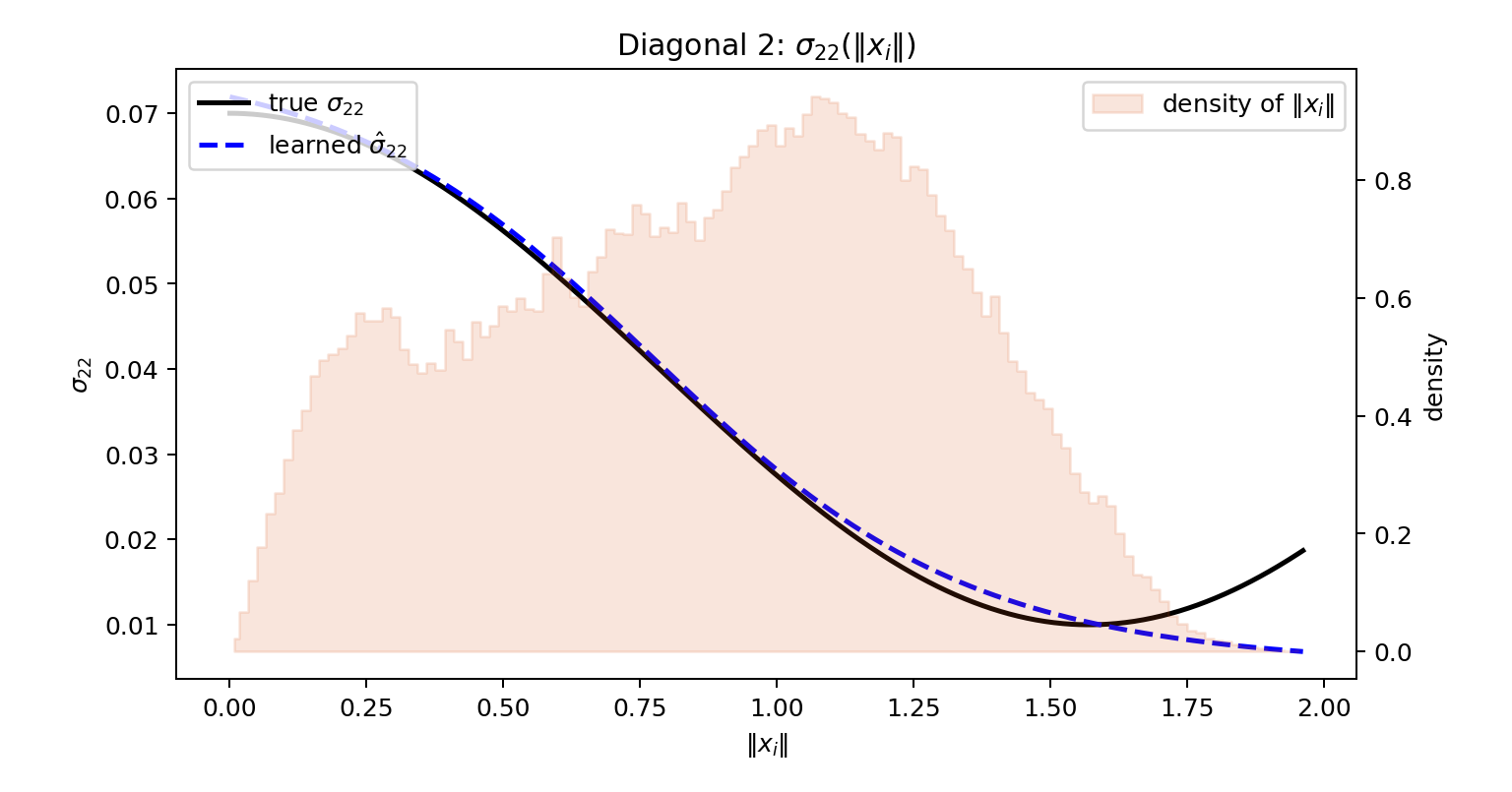}
  \caption{Case (II): Err = $0.025$.}
  \label{fig:Agent_sigma_2_22}
  \end{subfigure}
\caption{True \(\Sigma_{ii}\) vs learned \(\hat \Sigma_{ii}\) for $i = 1, 2$.}
\label{fig:Agent_sigma}
\end{figure}
\textbf{Conclusion}: the comparison of the trajectories in Fig.~\ref{fig:60D_trajs_1} and~\ref{fig:60D_trajs_2} shows that the learned $\hat\rvx$ is close to the true $\rvx$ under the same noise.  The comparison of $\phi$ vs $\hat\phi$ in Fig.~\ref{fig:60D_phi_1} and~\ref{fig:60D_phi_2} shows that when the data is abundant (the background shows the pairwise distance data used to obtain $\hat\phi$), the two are close to each other;  for $r$ close to zero, due to the form of the system, i.e. $\phi(\norm{\rvx_j - \rvx_i})(\rvx_j - \rvx_i)$, the information is weighted by zero, our learning is not that promising.  Figure \ref{fig:Agent_sigma} show our estimation result on state dependent \(\sigma\) under two different kinds of dynamics. Each diagonal entry is modeled by a shallow two–hidden–layer Tanh network with width $32$. The estimators tracks the true $\sigma$ closely even with such a lightweight network. 
\subsection{Example: SPDE estimation}\label{sec:spde_test}
We extend our method of section~\ref{sec:Learning Framework} to the stochastic heat equation with additive noise
\begin{equation}\label{eq:SPDE_general}
  \dif\rvu(t,\rvx)-\theta(\rvx)\,\Delta\rvu(t,\rvx)\,\dif t = \sigma\,\dif\rvw(t,\rvx),
\end{equation}
on a smooth bounded domain $\mathcal{D} \subset\R^{d}$, with initial condition $\rvu(0, \rvx) = 0$, zero boundary condition, and where $\Delta$ denotes the Laplace operator on $\mathcal{D} $ with zero boundary conditions. The existence, uniqueness and other analytical properties of the solution $\rvu$ are well understood, and we refer to~\citep{LototskyRozobsky2017Book}. Throughout this section, we fix the Hilbert space $H=L^{2}(\mathcal{D} )$ equipped with the usual inner product denoted by $(\cdot, \cdot)_H$. We note that in this case, the Laplace operator $\Delta$ has only point spectrum, and we denote by \( \bigl\{h_k:\,k\in\mathbb{N} \bigr\}\subset H\) its  eigenfunctions and $-\lambda_k$ the corresponding eigenvalues, i.e. \( \Delta h_k = -\lambda_k h_k \). It is well known that $\{h_k:k\in\mathbb{N}\}$ is a complete system in $H$, and without loss of generality we assume it is also orthonormal. The space-time noise, is assumed to be a cylindrical Brownian motion in $H$, which  informally can be written as
\(
  \rvw(t,\rvx)=\sum_{k\in\mathbb{N}} q_k h_k(\rvx)\,\rvw_k(t),
\)
where $\{q_k\}_{k\in\mathbb{N}}\subset(0,\infty)$ and $\{\rvw_k\}_{k\in\mathbb{N}}$ are independent one dimensional Brownian motions and $\sigma$ is a positive constant. Assume that $\theta$ is bounded, a.s. continuous on $\mathcal{D}$, and $\theta(x)\geq c_0>0$, for some positive real $c_0$. This guarantees the existence of the solution to \eqref{eq:SPDE_general} in an appropriate triple of Hilbert spaces. We are interested in the estimation of $\theta(x)$. 
\begin{figure}[ht!]
    \centering
    \includegraphics[width=0.9\linewidth]{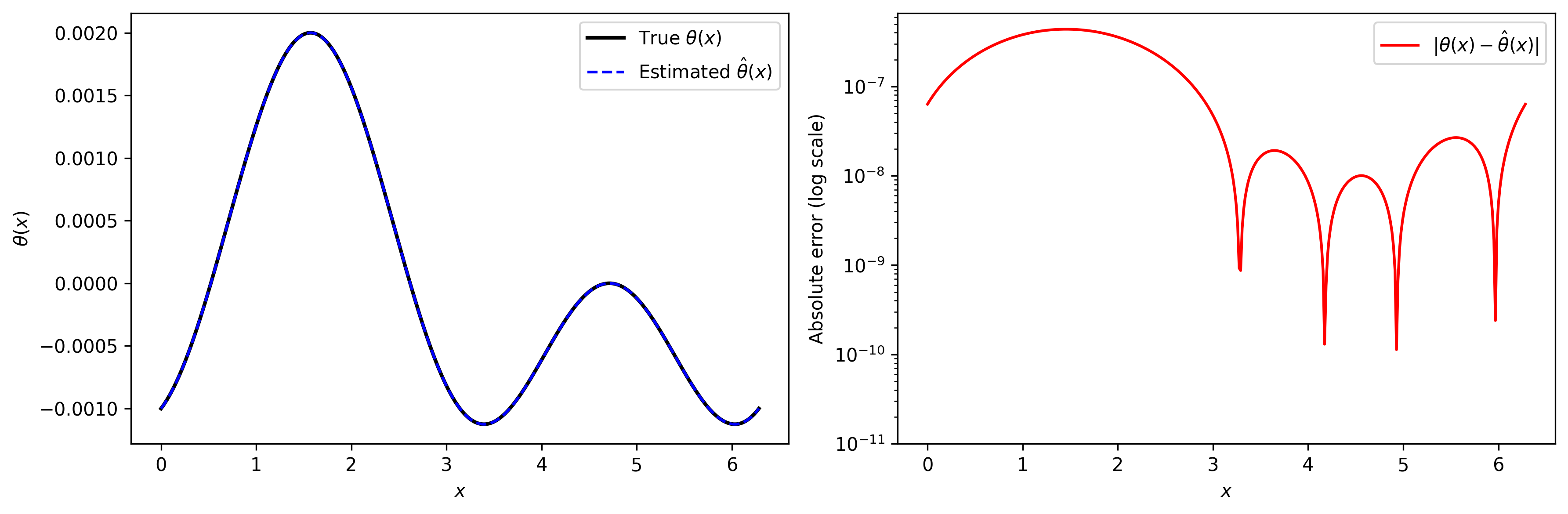}
    \caption{Left: Exact \(\theta_1(x)\) (solid) vs $\hat\theta_1$ (dashed). Right: \(|\theta_1 - \hat \theta_1|\) in log-scale.}
    \label{fig:SPDE_1}
\end{figure}
To verify our theoretical result, we present two numerical experiments for the stochastic heat equation~\eqref{eq:SPDE_general} on the spatial interval \([0,2\pi]\). Throughout we consider Fourier basis as our estimation function space, i.e., $\gH_{n}^\theta =\operatorname{span}\bigl\{1,\sin(k\rvx),\cos(k\rvx)\bigr\}_{k=1}^{n}$. In simulation of $\rvu(t, \rvx)$, we apply a Galerkin projection of dimension \(N\) with time step \(\Delta t = 10^{-3}\) up to horizon \(T = 10\). 
\begin{wrapfigure}{l}{0.5\textwidth}
\includegraphics[width=0.9\linewidth]{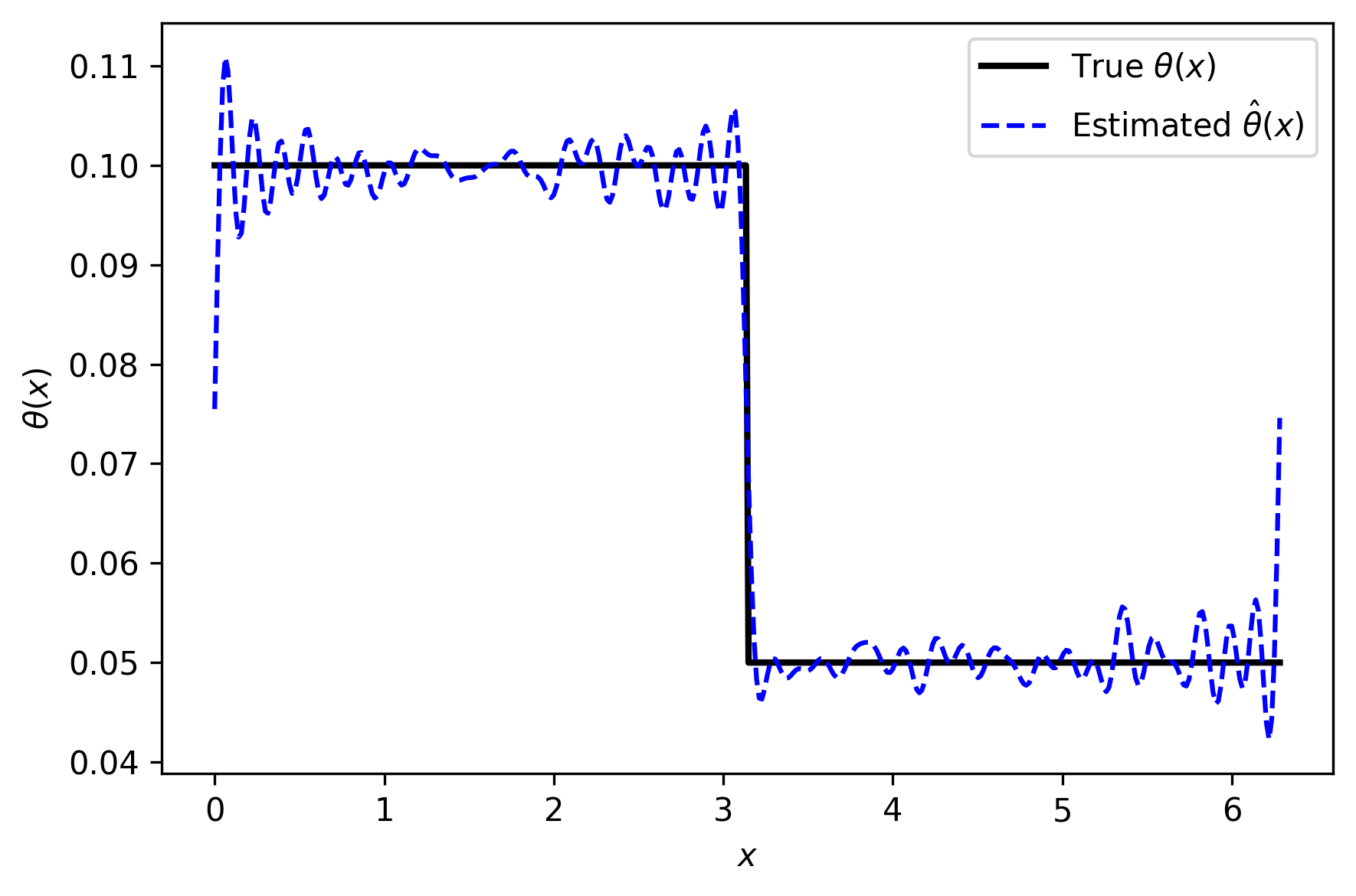} 
\caption{True \(\theta_2\) (solid) vs $\hat\theta_2$ (dashed).}
\label{fig:SPDE_2}
\end{wrapfigure}
The drift matrix is calculated with true \(\theta(x)\) where no projection error is introduced at the simulation stage.  Only the first \(N_{\mathrm{obs}}\) highest--frequency Fourier modes are marked as observable and the noise level is fixed at \(\sigma = 0.2\). In the first experiment, we take $\theta_1(x) \;=\; 0.001\,\bigl(\sin x - \cos 2x\bigr) \in \gH_2^\theta$, set \(N = 100\) and \(N_{\mathrm{o}} = 40\).   Next we consider a discontinuous coefficient outside the function space with
\[
\theta_2(x) \;=\;
\begin{cases}
        0.10, & 0 \leq x < \pi\\[3pt]
        0.05, & \pi \leq x \leq 2\pi
\end{cases}
\]
and set the estimation function subspace to \(\gH_{40}^\theta\). The simulation is carried out with \(N = 200\) and we observe \(N_{\mathrm{o}} = 100\) modes. Figure~\ref{fig:SPDE_1} and ~\ref{fig:SPDE_2} show the effectiveness of our learning under two fundamentally different scenarios, one with $\theta \in \gH_n^\theta$ and the other with $\theta \not\in \gH_n^\theta$.

\section{Additional Experiments}\label{sec:supp_examples}
In this section we report three $1D$ experiments and one $2D$ experiment designed to discuss $(i)$ the robustness to stochastic noise magnitude, $(ii)$ the effect of observation noise, $(iii)$ the effect of the learning time step (between the learning sampling instances and integration time instances), and $(iv)$ a correlated state-dependent noise structure where we show the effects of using a learned noise diffusion matrix as well as a comparison to traditional methods.  We also include a convergence study to empirically validate the statistical consistency predicted by our theory. For the $1D$ cases, we consider the following examples
\[
  \dif\rvx_t = \vf_\star(\rvx_t)\,\dif t + \sigma_\star(\rvx_t)\,\dif\rvw_t,
  \qquad \vf_\star(\rvx) = \rvx - \rvx^3,
\]
simulate trajectories by Euler-Maruyama with time step $\delta t = 10^{-3}$ on $[0,T]$ with $T=10$, and use $M=500$ independent trajectories. The drift is learned in the polynomial space using the discrete version of our noise-aware drift loss, and the error is measured in relative $L^2(\rho)$-norm.

\subsection{Experiment $1$: Varying the diffusion magnitude}
\begin{figure}[t]
  \centering
  \includegraphics[width=0.48\linewidth]{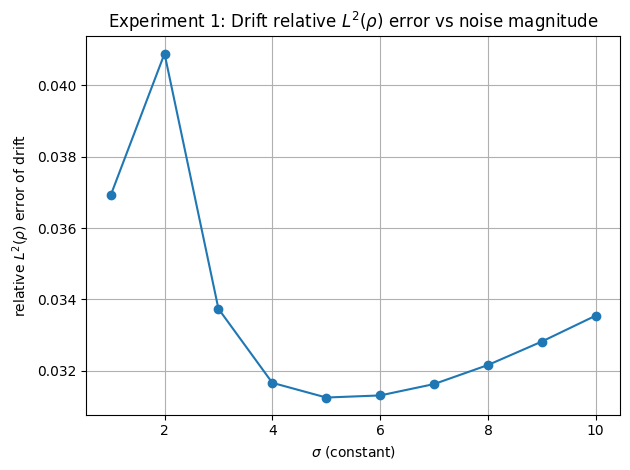}
  \includegraphics[width=0.48\linewidth]{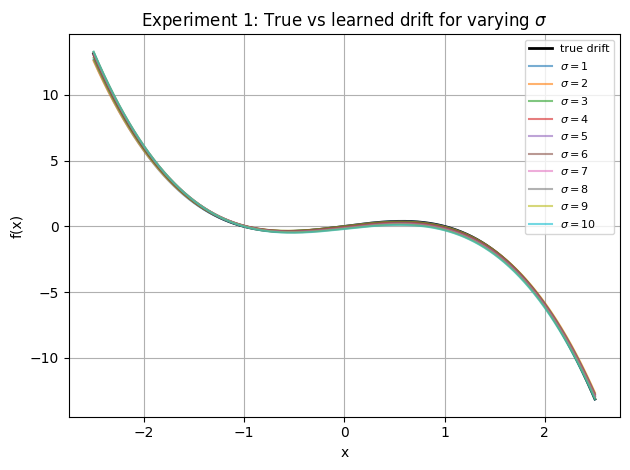}
  \caption{Experiment~1. Left: relative $L^2(\rho)$ error of the drift vs.\ constant diffusion level $\sigma$. Right: true drift $\vf_\star$ and learned drifts $\hat \vf_\sigma$ for several values of $\sigma$.}
  \label{fig:exp1}
\end{figure}

The goal of the first experiment is to test how the noise-aware drift estimator behaves as the dynamical noise level varies. We fix
\[
  \sigma_\star(x) = \sigma \in \{1,2,\dots,10\},
\]
simulate $\rvx_t$ with the true drift $\vf_\star$ and diffusion $\sigma_\star$. In the drift loss we treat $\sigma^2$ as known and plug in the true value. For each $\sigma$ we compute the learned drift $\hat \vf_\sigma$ and its relative $L^2(\rho)$ error. We report the results in Fig.~\ref{fig:exp1}, which suggests that the increasing $\sigma$ (stochastic noise level) does not affect our learning of the drift much.  In fact, even if $\sigma$ is a matrix, as long as its spectrum is bounded above zero, then our learning can handle it smoothly.

\subsection{Experiment $2$: Observation noise}
The second experiment investigates the effect of observation noise. We fix a smooth state-dependent diffusion,
\[
  \sigma_\star(\rvx) = 0.5 + 0.2\sin(\rvx),
\]
and simulate the true process $\rvx_t$ as before. On a learning grid with $\Delta t = 10^{-3}$ we form clean states $\rvx_{t_\ell}$ and then define noisy observations
\[
  \rvy_{t_\ell} = \rvx_{t_\ell} + \varepsilon_\ell,\qquad
  \varepsilon_\ell = \gamma\,\sigma_\star(\rvx_{t_\ell})\sqrt{\Delta t}\,Z_\ell,
  \quad Z_\ell \sim \mathcal{N}(0,1),
\]
so that $\Var(\varepsilon_\ell \mid \rvx_{t_\ell}) = \gamma^2 \sigma_\star^2(\rvx_{t_\ell}) \Delta t$ and the observation-noise contribution to increments $\Delta \rvy_{t_\ell}$ has variance of order $2\gamma^2\sigma_\star^2(\rvx_{t_\ell})\Delta t$. In the loss we learn the drift from the noisy increments $\Delta \rvy_{t_\ell} = \rvy_{t_{\ell+1}} - \rvy_{t_\ell}$, and vary $\gamma$ so that the observation-noise variance ranges from negligible to comparable with the SDE noise. 

\begin{figure}[t]
  \centering
  \includegraphics[width=0.48\linewidth]{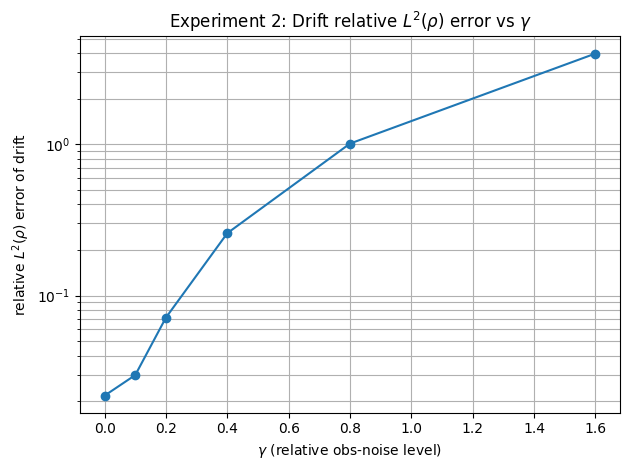}
  \includegraphics[width=0.48\linewidth]{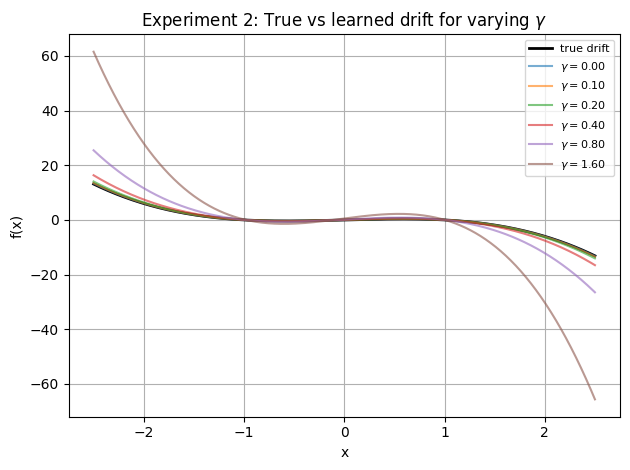}
  \caption{Experiment~2. Left: relative $L^2(\rho)$ drift error vs.\ observation-noise level $\gamma$. Right: true drift $\vf_\star$ and learned drifts $\hat \vf_\gamma$ for several values of $\gamma$.}
  \label{fig:exp2}
\end{figure}

This experiment illustrates that our estimator is robust when observation noise is small compared to the intrinsic SDE noise. And dealing with observation noise is filtering problem which is out of the scope of this paper. The results are shown in Fig.~\ref{fig:exp2}.

\subsection{Experiment $3$: Effect of the learning time step $\Delta t$}
The third experiment studies the impact of sampling step size. We fix a constant diffusion
\[
  \sigma_\star(x) \equiv 2,
\]
simulate with step size $\delta t = 10^{-3}$, and then subsample the trajectories on learning grids with
\[
  \Delta t \in \{0.1,\; 0.05,\; 0.025,\; 0.0125\}.
\]

\begin{figure}[t]
  \centering
  \includegraphics[width=0.48\linewidth]{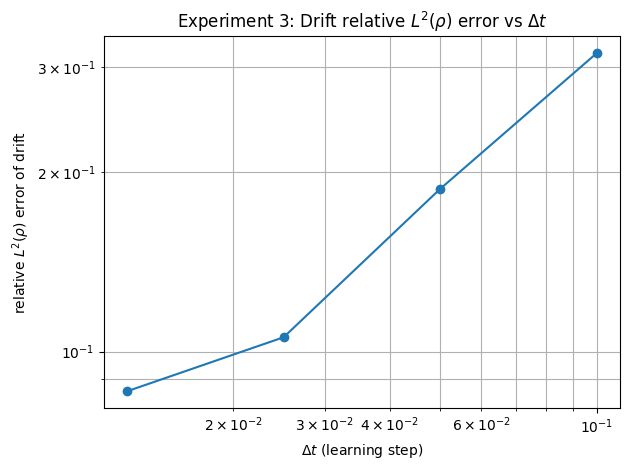}
  \includegraphics[width=0.48\linewidth]{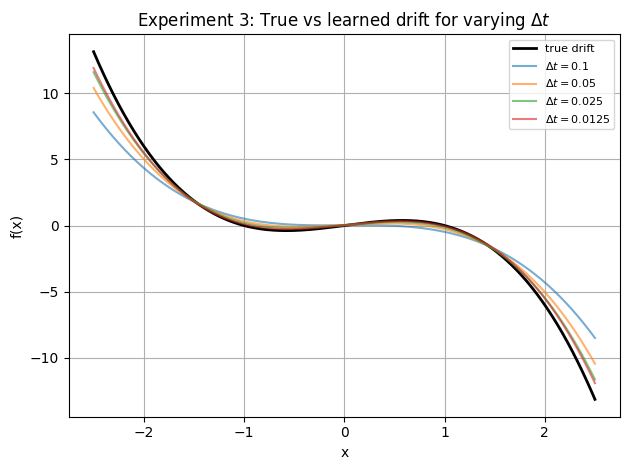}
  \caption{Experiment~3. Left: relative $L^2(\rho)$ drift error vs.\ learning step $\Delta t$. Right: True drift $\vf_\star$ and learned drifts $\hat \vf_{\Delta t}$ for several values of the learning step $\Delta t$.}
  \label{fig:exp3_errors}
\end{figure}
Together, these three 1D experiments presents the theoretical properties of our estimator: (i) robustness to changes in the intrinsic noise magnitude; (ii) the expected sensitivity to observation noise, which is not included in our research target; and (iii) different sampling size being consistent with the discrete-time approximation of the continuous-time loss. See Fig.~\ref{fig:exp3_errors} for the impact of the learning step size on the learned drift.

\subsection{Experiment $4$: $2$-Dim Stochastic Van der Pol with correlated state-dependent diffusion}\label{sec:methods_comp}
\begin{figure}[H]
  \centering
  \includegraphics[width=\linewidth]{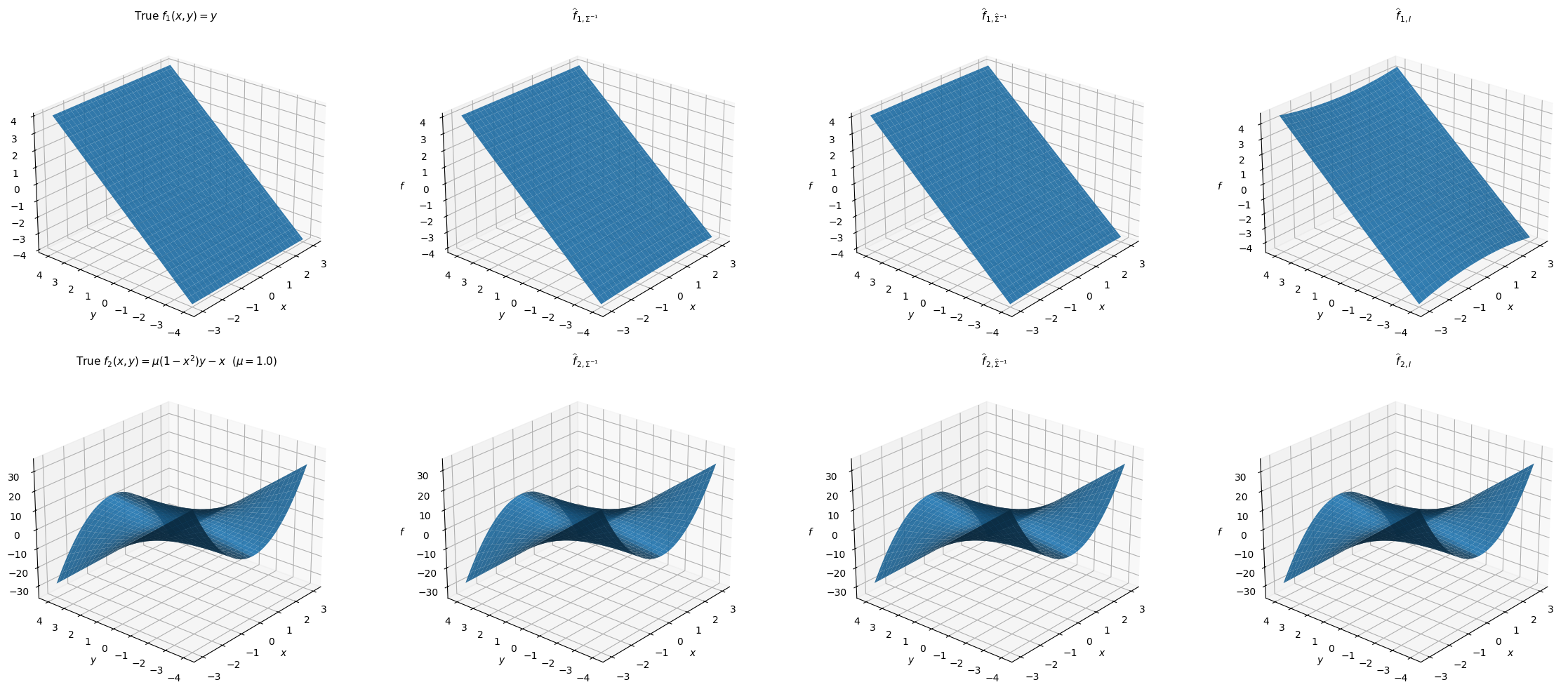}
  \caption{Experiment~4 (2D Van der Pol with correlated, state-dependent diffusion). 
  Top row: true $\vf_1(\rvx,\rvy)=\rvy$ and learned $\vf_1$ obtained with $\Sigma_\star^{-1}$, $\widehat\Sigma^{-1}$, and the SINDy-like baseline $\Sigma_\star=I$. 
  Bottom row: true $\vf_2(\rvx,\rvy)=\mu(1-\rvx^2)\rvy-\rvx$ and the corresponding learned $\vf_2$ for the same three choices.}
  \label{fig:exp4_vdp}
\end{figure}
We also include a 2D example based on the Van der Pol oscillator with a fully non-diagonal, state-dependent diffusion matrix. The drift is
\[
  \vf_\star(\rvx,\rvy) =
  \begin{pmatrix}
    \rvy \\
    \mu(1 - \rvx^2)\rvy - \rvx
  \end{pmatrix},
  \qquad \mu = 1,
\]
and we choose a volatility matrix $\sigma_\star(\rvx,\rvy) \in \mathbb{R}^{2\times 2}$ of the form
\[
  \Sigma_\star(\rvx,\rvy) \coloneqq \sigma_\star(\rvx,\rvy)\sigma_\star(\rvx,\rvy)^\top
  =
  \begin{pmatrix}
    v_1(\rvx)^2 & \rho(\rvx,\rvy)\,v_1(\rvx)\,v_2(\rvy)\\[0.2em]
    \rho(\rvx,\rvy)\,v_1(\rvx)\,v_2(\rvy) & v_2(\rvy)^2
  \end{pmatrix},
\]
with
\[
  v_1(\rvx) = 0.1 + 0.03\,\rvx^2,\qquad
  v_2(\rvy) = 0.2 + 0.04\,\rvy^2,\qquad
  \rho(\rvx,\rvy) = 0.5\,\tanh(0.2\,\rvx \rvy),
\]
so that $\Sigma_\star(\rvx,\rvy)$ is smooth, positive definite, and non-diagonal. We simulate $M=500$ trajectories on $[0,1]$ with time step $\delta t = 10^{-3}$.

We compare our noise aware learning with existing SINDy-like regression methods where $\Sigma_\star = I$ being an unweighted least-squares loss, see section \ref{sec:Comparison} for details. This is exactly the structure used in traditional SINDy-type drift estimators, which ignore the correlated, state-dependent covariance. The comparison is summarized in Fig.~\ref{fig:exp4_vdp}.

\subsection{Convergence Study}
\begin{figure}[ht!]
  \begin{subfigure}[b]{0.45\textwidth}
  \centering
  \includegraphics[width=\textwidth]{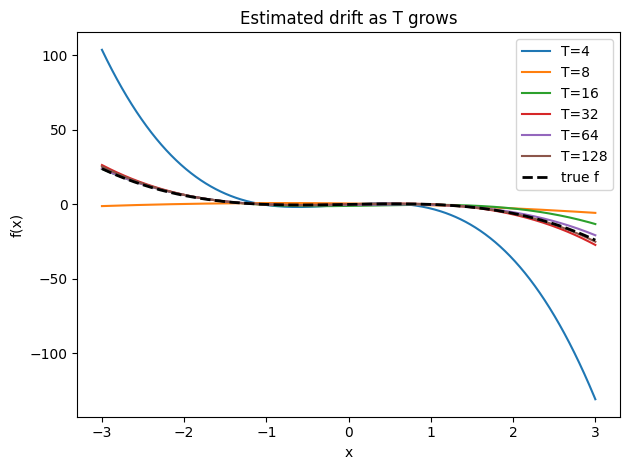}
  \caption{Estimated drifts (colored) vs true drift (black).}
  \label{fig:convergence test of f on T}
  \end{subfigure}
  \hfill  
  \begin{subfigure}[b]{0.45\textwidth}
  \centering
  \includegraphics[width=\textwidth]{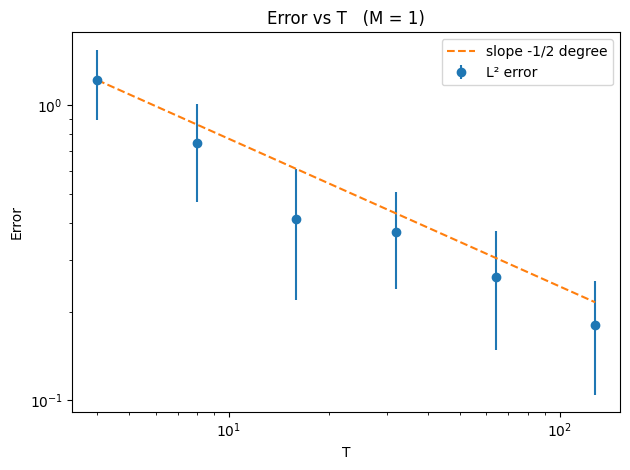}
  \caption{Log–log of $\|\hat \vf-\vf\|_{L^{2}(\rho)}$ vs $T$.}
  \label{fig:convergence test error on T}
  \end{subfigure}
\caption{Convergence Test with $M = 1$.}
\label{fig:T_conv}
\end{figure}

\begin{figure}[ht!]
  \begin{subfigure}[b]{0.45\textwidth}
  \centering
  \includegraphics[width=\textwidth]{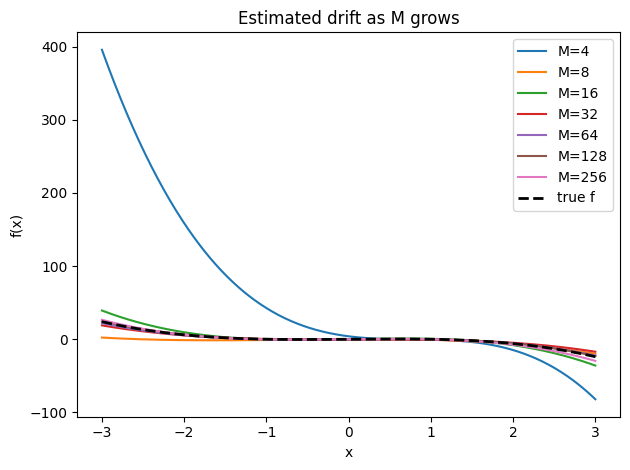}
  \caption{Estimated drifts (colored) vs true drift (black).}
  \label{fig:convergence test of f on M}
  \end{subfigure}
\hfill  
  \begin{subfigure}[b]{0.45\textwidth}
  \centering
  \includegraphics[width=\textwidth]{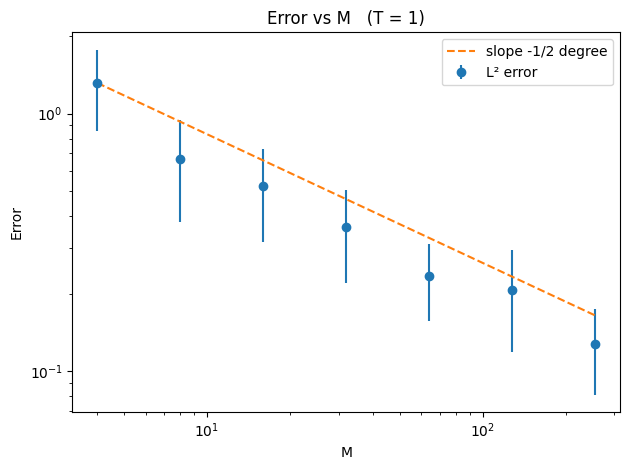}
  \caption{Log–log of $\|\hat \vf-\vf\|_{L^{2}(\rho)}$ vs $M$.}
  \label{fig:convergence test error on M}
  \end{subfigure}
\caption{Convergence Test with $T = 1$.}
\label{fig:M_conv}
\end{figure}

To illustrate the statistical consistency of our estimator defined in \ref{sec:Learning Framework}, we
consider an SDE in $d = 1$ case with
\(\vf(\rvx)=-\rvx^{3}+\rvx\) and \(\sigma(\rvx) = 1+0.4\,\sin \rvx\) simulated by the Euler–Maruyama scheme with step--size $\Delta t=10^{-3}$.
The initial states are drawn i.i.d.\ from the invariant density, so the process is strictly stationary. In $1D$, this density is
\[
\pi(\rvx)=\frac{1}{G(\vf)}\,\sigma(\rvx)^{-2}\exp\!\Big\{\,2\!\int_{0}^{\rvx}\frac{\vf(\rvv)}{\sigma(\rvv)^{2}}\,\dif \rvv\Big\}, \quad
G(\vf)=\int_{\R}\sigma(\rvx)^{-2}\exp\!\Big\{\,2\!\int_{0}^{\rvx}\frac{\vf(\rvv)}{\sigma(\rvv)^{2}}\,\dif \rvv\Big\}\,\dif \rvx .
\] 
For a collection of $M$ paths observed over $[0,T]$ our drift estimator $\hat{\vf}$ is searched in the space $\gH =  \operatorname{span}\{1,\rvx,\rvx^{2},\rvx^{3}\}$ by minimizing the loss function \ref{eq:original_loss}. The estimation error is quantified in the $\rho$–weighted norm introduced in
\eqref{eq:rho_norm}. We test consistency in both time $T$ and number of observed trajectories $M$ with each replicated $20$ times to obtain error bars. We first fix $M=1$ and let $T\in\{4,8,16,32,64,128\}$. Due to the ergodicity of the underlying SDE, we expect the following convergence rate $\lVert\hat \vf-\vf\rVert_{L^{2}(\rho)}=O(T^{-1/2})$, which is confirmed by Fig.\ref{fig:convergence test error on T}.
Next, we fix $T=1$ and let $M\in\{4,8,16,32,64,128,256\}$. The error decays at the rate $\lVert\hat \vf-\vf\rVert_{L^{2}(\rho)}=O(M^{-1/2})$, which is the rate confirmed by our theorem; see Fig.\ref{fig:convergence test error on M}.
In addition to Log–log plots~\ref{fig:convergence test error on T} and~\ref{fig:convergence test error on M} confirming  the predicted slopes $-1/2$ in both regimes, we plot the corresponding drift functions~\ref{fig:convergence test of f on T} and \ref{fig:convergence test of f on M} to illustrate the qualitative tightening of $\hat \vf$ towards $\vf$ as information increases.
These numerical findings demonstrate
that the estimator remains statistically consistent when the
diffusion coefficient is state dependent.
\section{Conclusion}\label{sec:conclude}
We have demonstrated a novel learning methodology for inferring the drift and diffusion coefficient in general SDE systems driven by Brownian noise.  Our estimation approach does not assume a specific functional structure for the drift or the diffusion, thereby enhancing its applicability across a diverse range of SDE and SPDE models. This approach can handle high-dimensional SDE systems by leveraging deep learning architectures. The loss function for the drift is derived from the negative logarithm of the ratio of likelihood functions. For the diffusion coefficient, the loss function is based on the quadratic variation, which operates independently of the drift function. This independence makes our method particularly effective in scenarios where only trajectory observations are available. Additionally, our approach is adaptable to various noise structures.  

\subsubsection*{Author Contributions}
ZG developed the algorithm, analyzed the data and implemented the software package.  ZG develops the theores along with IC and MZ.  IC and MZ designed the research.  All authors wrote the manuscript. 
\subsubsection*{Acknowledgments}
IC research was partially supported by NSF Grant DMS-2407549. MZ gratefully acknowledges funding provided by the Oak Ridge Associated Universities (ORAU) Ralph E. Powe Junior Faculty Enhancement Award and NSF Grant CCF-AF-$2225507$.  
\bibliography{SDE_Esti}

@article{HER2024,
title = {Pseudo-Hamiltonian system identification},
journal = {Journal of Computational Dynamics},
volume = {11},
number = {1},
pages = {59-91},
year = {2024},
issn = {2158-2491},
doi = {10.3934/jcd.2024001},
author = {Sigurd Holmsen and Sølve Eidnes and Signe Riemer-Sørensen},
}

@article{MIS2023,
title = {Learning dynamical systems using local stability priors},
journal = {Journal of Computational Dynamics},
volume = {10},
number = {1},
pages = {175-198},
year = {2023},
issn = {2158-2491},
doi = {10.3934/jcd.2022021},
author = {Arash Mehrjou and Andrea Iannelli and Bernhard Schölkopf},
}

@book{Sarkka_Solin_2019,
  place     = {Cambridge},
  series    = {Institute of Mathematical Statistics Textbooks},
  title     = {Applied Stochastic Differential Equations},
  publisher = {Cambridge University Press},
  author    = {Simo S\"{a}rkk\"{a} and Arno Solin},
  year      = {2019},
  collection= {Institute of Mathematical Statistics Textbooks}
}

@book{Evans,
  title     = {An Introduction to Stochastic Differential Equations},
  publisher = {American Mathematical Society},
  author    = {Lawrence C. Evans},
  year      = {2013}
}

@ARTICLE{lu2022,
  author  = {Fei Lu and Mauro Maggioni and Sui Tang},
  title   = {Learning interaction kernels in stochastic systems of interacting particles from multiple trajectories},
  journal = {Foundations of Computational Mathematics},
  year    = {2022},
  volume  = {22},
  pages   = {1013 -- 1067}
}

@ARTICLE{EGNS2008,
  author  = {Werner Ebeling and Ewa Gudowska-Nowak and Igor M. Sokolov},
  title   = {On Stochastic Dynamics in Physics - Remarks on History and Terminology},
  journal = {Acta Physica Polonica B},
  year    = {2008},
  volume  = {39},
  number  = {5},
  pages   = {1003 -- 1018}
}

@article{wanner2024higher,
  title={On higher order drift and diffusion estimates for stochastic SINDy},
  author={Wanner, Mathias and Mezi{\'c}, Igor},
  journal={SIAM Journal on Applied Dynamical Systems},
  volume={23},
  number={2},
  pages={1504--1539},
  year={2024},
  publisher={SIAM}
}

@inproceedings{wen2022transformers,
  title     = {Transformers in Time Series: A Survey},
  author    = {Qingsong Wen and Tian Zhou and Chaoli Zhang and Weiqi Chen and Ziqing Ma and Junchi Yan and Liang Sun},
  booktitle = {Proceedings of the Thirty-Second International Joint Conference on Artificial Intelligence (IJCAI 2023), Survey Track},
  year      = {2023},
  doi       = {10.24963/ijcai.2023/759}
}

@article{ho2020denoising,
  title   = {Denoising diffusion probabilistic models},
  author  = {Jonathan Ho and Ajay Jain and Pieter Abbeel},
  journal = {Advances in Neural Information Processing Systems},
  volume  = {33},
  pages   = {6840--6851},
  year    = {2020}
}

@inproceedings{song2020score,
  title     = {Score-based generative modeling through stochastic differential equations},
  author    = {Yang Song and Jascha Sohl-Dickstein and Diederik P. Kingma and Abhishek Kumar and Stefano Ermon and Ben Poole},
  booktitle = {Proceedings of the 9th International Conference on Learning Representations (ICLR)},
  year      = {2021}
}

@ARTICLE{DP2011,
  author  = {David Dingli and Jorge M. Pacheco},
  title   = {Stochastic dynamics and the evolution of mutations in stem cells},
  journal = {BMC Biol},
  year    = {2011},
  volume  = {9},
  pages   = {41}
}

@ARTICLE{sindy,
  author  = {Steven L. Brunton and Joshua L. Proctor and J. Nathan Kutz},
  title   = {Discovering governing equations from data by sparse identification of nonlinear dynamical systems},
  journal = {PNAS},
  year    = {2016},
  volume  = {113},
  number  = {15},
  pages   = {3932 -- 3937}
}

@article{chen2018neural,
  title   = {Neural ordinary differential equations},
  author  = {Ricky T. Q. Chen and Yulia Rubanova and Jesse Bettencourt and David K. Duvenaud},
  journal = {Advances in Neural Information Processing Systems},
  volume  = {31},
  year    = {2018}
}

@article{pinn2019,
  title     = {Physics-informed neural networks: A deep learning framework for solving forward and inverse problems involving nonlinear partial differential equations},
  journal   = {Journal of Computational Physics},
  volume    = {378},
  pages     = {686-707},
  year      = {2019},
  issn      = {0021-9991},
  doi       = {https://doi.org/10.1016/j.jcp.2018.10.045},
  url       = {https://www.sciencedirect.com/science/article/pii/S0021999118307125},
  author    = {Maziar Raissi and Paris Perdikaris and George E. Karniadakis}
}

@article{CS02,
  title   = {On the mathematical foundations of learning},
  author  = {Felipe Cucker and Steve Smale},
  journal = {Bulletin of the American Mathematical Society},
  volume  = {39},
  number  = {1},
  pages   = {1-49},
  year    = {2002}
}

@article{0b9b8115-a8b8-3422-8e1c-a62077de6621,
  ISSN     = {00223808, 1537534X},
  URL      = {http://www.jstor.org/stable/1831029},
  author   = {Fischer Black and Myron Scholes},
  journal  = {Journal of Political Economy},
  number   = {3},
  pages    = {637--654},
  publisher= {University of Chicago Press},
  title    = {The Pricing of Options and Corporate Liabilities},
  urldate  = {2023-12-11},
  volume   = {81},
  year     = {1973}
}

@article{VASICEK1977177,
  title   = {An equilibrium characterization of the term structure},
  journal = {Journal of Financial Economics},
  volume  = {5},
  number  = {2},
  pages   = {177-188},
  year    = {1977},
  issn    = {0304-405X},
  doi     = {https://doi.org/10.1016/0304-405X(77)90016-2},
  url     = {https://www.sciencedirect.com/science/article/pii/0304405X77900162},
  author  = {Oldrich Vasicek}
}

@article{10.1093/rfs/6.2.327,
  author  = {Steven L. Heston},
  title   = {A Closed-Form Solution for Options with Stochastic Volatility with Applications to Bond and Currency Options},
  journal = {The Review of Financial Studies},
  volume  = {6},
  number  = {2},
  pages   = {327-343},
  year    = {1993},
  month   = {04},
  issn    = {0893-9454},
  doi     = {10.1093/rfs/6.2.327},
  url     = {https://doi.org/10.1093/rfs/6.2.327},
  eprint  = {https://academic.oup.com/rfs/article-pdf/6/2/327/24417457/060327.pdf}
}

@ARTICLE{935092,
  author  = {Yaser S. Abu-Mostafa},
  journal = {IEEE Transactions on Neural Networks},
  title   = {Financial model calibration using consistency hints},
  year    = {2001},
  volume  = {12},
  number  = {4},
  pages   = {791-808},
  doi     = {10.1109/72.935092}
}

@book{doi:10.1142/8195,
  author    = {William T. Coffey and Yuri P. Kalmykov},
  title     = {The Langevin Equation},
  publisher = {WORLD SCIENTIFIC},
  year      = {2012},
  doi       = {10.1142/8195},
  edition   = {3rd},
  URL       = {https://www.worldscientific.com/doi/abs/10.1142/8195},
  eprint    = {https://www.worldscientific.com/doi/pdf/10.1142/8195}
}

@article{SZEKELY201414,
  title    = {Stochastic simulation in systems biology},
  journal  = {Computational and Structural Biotechnology Journal},
  volume   = {12},
  number   = {20},
  pages    = {14-25},
  year     = {2014},
  issn     = {2001-0370},
  doi      = {https://doi.org/10.1016/j.csbj.2014.10.003},
  url      = {https://www.sciencedirect.com/science/article/pii/S2001037014000403},
  author   = {Tam{\'a}s Sz{\'e}kely and Kevin Burrage}
}

@article{10.1063/1.4948407,
  author  = {Fuke Wu and Tianhai Tian and James B. Rawlings and George Yin},
  title   = {Approximate method for stochastic chemical kinetics with two-time scales by chemical {L}angevin equations},
  journal = {The Journal of Chemical Physics},
  volume  = {144},
  number  = {17},
  pages   = {174112},
  year    = {2016},
  month   = {05},
  issn    = {0021-9606},
  doi     = {10.1063/1.4948407},
  url     = {https://doi.org/10.1063/1.4948407},
  eprint  = {https://pubs.aip.org/aip/jcp/article-pdf/doi/10.1063/1.4948407/13987287/174112_1_online.pdf}
}

@article{sachs2017langevin,
  title     = {Langevin dynamics with variable coefficients and nonconservative forces: from stationary states to numerical methods},
  author    = {Matthias Sachs and Benedict Leimkuhler and Vincent Danos},
  journal   = {Entropy},
  volume    = {19},
  number    = {12},
  pages     = {647},
  year      = {2017},
  publisher = {MDPI}
}

@article{takeuchi2006evolution,
  title    = {Evolution of predator--prey systems described by a {L}otka--{V}olterra equation under random environment},
  author   = {Yasuhiro Takeuchi and Nguyen Huy Du and Nguyen Thai Hieu and Kazunori Sato},
  journal  = {Journal of Mathematical Analysis and Applications},
  volume   = {323},
  number   = {2},
  pages    = {938--957},
  year     = {2006},
  publisher= {Elsevier}
}

@article{liao2019learning,
  title   = {Learning stochastic differential equations using {RNN} with log signature features},
  author  = {Shujian Liao and Terry Lyons and Weixin Yang and Hao Ni},
  journal = {arXiv preprint arXiv:1908.08286},
  year    = {2019}
}

@article{YANG2023279,
  title   = {Neural network stochastic differential equation: models with applications to financial data forecasting},
  journal = {Applied Mathematical Modelling},
  volume  = {115},
  pages   = {279-299},
  year    = {2023},
  issn    = {0307-904X},
  doi     = {https://doi.org/10.1016/j.apm.2022.11.001},
  url     = {https://www.sciencedirect.com/science/article/pii/S0307904X22005340},
  author  = {Luxuan Yang and Ting Gao and Yubin Lu and Jinqiao Duan and Tao Liu}
}

@book{hull2017options,
  title     = {Options, Futures, and Other Derivatives},
  author    = {John C. Hull},
  isbn      = {9780134631493},
  url       = {https://books.google.com/books?id=yfr0DQAAQBAJ},
  year      = {2017},
  publisher = {Pearson Education}
}

@inproceedings{guo2024learning,
  title     = {Learning Stochastic Dynamics from Data},
  author    = {Ziheng Guo and Ming Zhong and Igor Cialenco},
  booktitle = {ICLR 2024 Workshop on AI4DifferentialEquations In Science},
  year      = {2024},
  url       = {https://openreview.net/forum?id=MdXtFDhy0H}
}

@article{doi:10.1137/S0036144500378302,
  author  = {Desmond J. Higham},
  title   = {An Algorithmic Introduction to Numerical Simulation of Stochastic Differential Equations},
  journal = {SIAM Review},
  volume  = {43},
  number  = {3},
  pages   = {525-546},
  year    = {2001},
  doi     = {10.1137/S0036144500378302},
  URL     = {https://doi.org/10.1137/S0036144500378302},
  eprint  = {https://doi.org/10.1137/S0036144500378302}
}

@article{LZTM2019,
  author  = {Fei Lu and Ming Zhong and Sui Tang and Mauro Maggioni},
  title   = {Nonparametric inference of interaction laws in systems of agents from trajectory data},
  year    = {2019},
  volume  = {116},
  number  = {29},
  pages   = {14424 -- 14433},
  journal = {Proceedings of the National Academy of Sciences}
}

@article{ZMM2020,
  author  = {Ming Zhong and Jason Miller and Mauro Maggioni},
  title   = {Data-driven Discovery of Emergent Behaviors in Collective Dynamics},
  year    = {2020},
  volume  = {411},
  month   = {October},
  pages   = {132542},
  journal = {Physica D: nonlinear phenomenon}
}

@article{FMMZ2022,
  title   = {Learning Interaction Variables and Kernels from Observations of Agent-Based Systems},
  journal = {IFAC-PapersOnLine},
  volume  = {55},
  number  = {30},
  pages   = {162-167},
  year    = {2022},
  note    = {25th International Symposium on Mathematical Theory of Networks and Systems MTNS 2022},
  author  = {Jinchao Feng and Mauro Maggioni and Patrick Martin and Ming Zhong}
}

@incollection{FM2023,
  author    = {Jinchao Feng and Ming Zhong},
  title     = {Learning Collective Behaviors from Observation},
  pages     = {},
  editor    = {Simon Foucart and Stephan Wojtowytsch},
  booktitle = {Explorations in the Mathematics of Data Science},
  booksubtitle = {The Inaugural Volume of the Center for Approximation and Mathematical Data Analytics},
  year      = {2024},
  publisher = {Birkh{"a}user Cham}
}

@InProceedings{MMQZ2021,
  title     = {Learning Interaction Kernels for Agent Systems on Riemannian Manifolds},
  author    = {Mauro Maggioni and Jason D. Miller and Hongda. Qiu and Ming Zhong},
  booktitle = {Proceedings of the 38th International Conference on Machine Learning},
  pages     = {7290--7300},
  year      = {2021},
  editor    = {Meila, Marina and Zhang, Tong},
  volume    = {139},
  series    = {Proceedings of Machine Learning Research},
  month     = {18--24 Jul},
  publisher = {PMLR}
}

@article{xu2023modeling,
  title   = {Modeling Unknown Stochastic Dynamical System via Autoencoder},
  author  = {Zhongshu Xu and Yuan Chen and Qifan Chen and Dongbin Xiu},
  journal = {Journal of Machine Learning for Modeling and Computing},
  year    = {2024},
  volume  = {5},
  number  = {3},
  pages   = {87--112}
}

@article{Talay2002StochasticHS,
  author  = {Daniel Talay},
  year    = {2002},
  title   = {Stochastic {H}amiltonian Systems: Exponential Convergence to the Invariant Measure, and Discretization by the Implicit {E}uler Scheme},
  volume  = {8},
  journal = {Markov Processes and Related Fields}
}

@book{liptser2001statistics,
  title     = {Statistics of Random Processes: I. General Theory},
  author    = {Robert S. Liptser and Albert N. Shiryaev},
  isbn      = {9783540639299},
  lccn      = {00041918},
  series    = {Applications of Mathematics Stochastic Modelling and Applied Probability Series},
  url       = {https://books.google.com/books?id=gKtK0CjxOaIC},
  year      = {2001},
  publisher = {Springer}
}

@Book{KutoyantsBook2004,
  author    = {Yury A. Kutoyants},
  publisher = {Springer-Verlag London Ltd.},
  title     = {Statistical Inference for Ergodic Diffusion Processes},
  year      = {2004},
  address   = {London},
  series    = {Springer Series in Statistics}
}

@book{LototskyRozobsky2017Book,
  title     = {Stochastic Partial Differential Equations},
  ISBN      = {9783319586472},
  ISSN      = {2191-6675},
  url       = {http://dx.doi.org/10.1007/978-3-319-58647-2},
  DOI       = {10.1007/978-3-319-58647-2},
  journal   = {Universitext},
  publisher = {Springer International Publishing},
  author    = {Sergey V. Lototsky and Boris L. Rozovsky},
  year      = {2017}
}

@article{MrazekPospisil2017,
  author  = {Milan Mr\'azek and Jan Posp\'i\v{s}il},
  title   = {Calibration and Simulation of {H}eston Model},
  journal = {Open Mathematics},
  year    = {2017},
  volume  = {15},
  number  = {1},
  pages   = {679--704},
  doi     = {10.1515/math-2017-0058}
}

@book{Vaart_1998, 
place={Cambridge}, 
series={Cambridge Series in Statistical and Probabilistic Mathematics}, 
title={Asymptotic Statistics},
publisher={Cambridge University Press}, 
author={Vaart, A. W. van der}, 
year={1998}, 
collection={Cambridge Series in Statistical and Probabilistic Mathematics}}

@incollection{NEWEY19942111,
title = {Chapter 36 Large sample estimation and hypothesis testing},
series = {Handbook of Econometrics},
publisher = {Elsevier},
volume = {4},
pages = {2111-2245},
year = {1994},
issn = {1573-4412},
doi = {https://doi.org/10.1016/S1573-4412(05)80005-4},
url = {https://www.sciencedirect.com/science/article/pii/S1573441205800054},
author = {Whitney K. Newey and Daniel McFadden}
}
\bibliographystyle{iclr2026_conference}
\appendix
\section{Learning Framework}
We discussion additional details related to the learning of drift and noise in this section.
\subsection{Simplification of the Loss}
When $D \gg 1$ and $\sigma = \sigma(\rvx) \in \R^{D \times D}$ is a full matrix, the learning of the drift term $\vf$ can be computationally demanding, as all components of $\vf$ are coupled and one has to solve the optimization problem in high-dimensional space all at once.  Stochastic gradient descent coupled with neural network solutions is one of the desired approaches; however the solutions become less interpretable.  In this section, we discuss several scenarios this loss for learning drift can be simplified.  In this section, we discuss several scenarios in which the loss for learning drift can be simplified.

In the case of the noise being a constant full matrix, i.e. $\sigma(\rvx_t) = \sigma \in \R^{D\times D}$, the loss is equivalent (in the optimization sense) to the following 
\[
\loss_{\gH}^{\text{Sim}}(\tilde\vf) = \E\Big[\int_{t = 0}^T\norm{\tilde\vf(\rvx_t)}^2\,\dt - 2\langle\tilde\vf(\rvx_t), \dif \rvx_t\rangle\Big]
\]

In the case of state-dependent uncorrelated noise, i.e. $\Sigma(\vx) = \sigma^2(\vx)\mI$, where $\mI$ is the $D \times D$ identity matrix and $\sigma: \R^D \rightarrow \R^+$ is a scalar function depending on the state and representing the noise level, the loss function \eqref{eq:original_loss} can be simplified to
\begin{equation}\label{eq:simple_loss}
\loss_{\gH}^{\text{Sim}}(\tilde\vf) =\E\Big[\sum_{d = 1}^D\int_{t = 0}^T\frac{|\tilde{f}_d(\rvx_t)|^2\,\dt - 2\tilde{f}_d(\dif\rvx)_d(t)}{2\sigma^2(\rvx_t)}\Big],
\end{equation}
where $\tilde\vf(\rvx_t) = (\tilde{f}_1(\rvx_t), \cdots, \tilde{f}_D(\rvx_t))$.  Hence the learning of each component of $\vf$ can be de-coupled.  When $\Sigma$ is a state-dependent full matrix, we consider the eigen-decomposition of $\Sigma$, i.e. $\Sigma(\vx) = \mQ\mLambda(\vx)\mQ^\intercal$, then we rotate the system by $\mQ^\intercal$, i.e., $\rvx_t' = \mQ^\intercal\rvx_t$, $\vf'(\vx') = \mQ^\intercal\vf(\vx)\mQ$, $\rvw_t' = \mQ^\intercal\rvw_t$, then we obtain the case when $\Sigma$ is a diagonal matrix.  Once we learn $\hat\mLambda$ and $\vf'$, we will use the following to obtain the original functions, i.e., $\vf(\vx) = \mQ\vf'(\mQ\vx)\mQ^\intercal$, and $\hat\Sigma = \mQ\hat\mLambda\mQ^\intercal$.
\subsection{Implementation}\label{subsec: Implementation}
We discuss in details how the algorithm is implemented for our learning framework. Practically speaking, data are rarely sampled continuously in time. Instead, observers typically have access to fragmented data sets, gathered from multiple independently sampled trajectories at specific, discrete time points$\{\rvx_l^m\}_{l, m = 1}^{L, M}$, where $\rvx_l^m = \rvx^{(m)}(t_l)$ with $0 = t_1 < \cdots < t_L = T$ and $\rvx_0^m$ is an i.i.d sample from $\mu_0$.  
We use a discretized version of \ref{eq:original_loss}, 
\begin{equation}\label{eq:original_loss_dis}
\loss_{L, M, \gH}(\tilde\vf) = \frac{1}{2M}\sum_{l, m = 1}^{L - 1, M}
\Big(\langle\tilde\vf(\rvx_l^m), \Sigma^{-1}(\rvx_l^m)\tilde\vf(\rvx_l^m)\rangle\Delta t_l  - 2\langle\tilde\vf(\rvx_l^m), \Sigma^{-1}(\rvx_l^m)\Delta\rvx_l^m\rangle\Big),
\end{equation}
for $\tilde\vf \in \gH$ and $\Delta\rvx_l^m = \rvx_{l + 1}^m - \rvx_l^m$.  Moreover, we also assume that $\gH$ is a finite-dimensional function space, i.e. $\Dim(\gH) = n < \infty$.  Then for any $\tilde\vf \in \gH$, $\tilde\vf(\vx) = \sum_{i = 1}^n\va_{i}\psi_{i}(\vx)$, where $\va_{i} \in \R^D$ is a constant vector coefficient and $\psi_{i}: \mD \subset \R^D \rightarrow \R$ is a basis of $\gH$ and the domain $\mD$ is constructed by finding out the $\min/\max$ of the components of $\rvx_t \in \R^D$ for $t \in [0, T]$.  We consider two methods for constructing $\psi_i$: $a)$ use pre-determined basis such as piecewise polynomials or Clamped B-spline, Fourier basis, or a mixture of all of the aforementioned ones; $b)$  use neural networks, where the basis functions are also trained from data.  Next, we can put the basis representation of $\tilde\vf$ back to \eqref{eq:original_loss_dis}, we obtain the following loss based on the coefficients
\begin{equation}\label{eq:original_loss_dis_coeff}
\loss_{L, M, \gH}(\{\va_{\eta}\}_{i = 1}^n) = \frac{1}{2M}\sum_{l, m = 1}^{L - 1, M}\Big(\sum_{i, j = 1}^n\langle\va_i, \Sigma^{-1}_{l, m}\va_j\rangle\psi^m_{i, l}\psi^m_{j, l}\Delta t_{l} - 2\sum_{i = 1}^n\langle\va_i, \Sigma^{-1}_{l, m}\Delta\rvx_l^m\rangle\psi^m_{i, l}\Big),
\end{equation}
where $\psi^m_{i, l} = \psi_i(\rvx_l^m)$, $\Sigma^{-1}_{l, m} = \Sigma^{-1}(\rvx_l^m)$ and $\Delta t_l = t_{l + 1} - t_l$. In the case of diagonal covariance matrix $\Sigma$, i.e., $\Sigma(\vx) = \text{diag}(\sigma_1^2(\vx), \cdots, \sigma_D^2(\vx)) \in \R^{D \times D}$, for $\sigma_i > 0$ and $i = 1, \cdots, D$;  we can re-write \eqref{eq:original_loss_dis_coeff}  as 
\[
\loss_{L, M, \gH}(\{\va_{\eta}\}_{i = 1}^n)= \frac{1}{2M}\sum_{l, m = 1}^{L - 1, M}\Big(\sum_{i, j}^{n}\frac{\langle\va_i, \va_j\rangle}{\sigma_k^2(\rvx_l^m)}\psi^m_{i, l}\psi^m_{j, l}\Delta t_l - 2\sum_{i = 1}^n\frac{\langle\va_i, \Delta\rvx_l^m\rangle}{\sigma_k^2(\rvx_l^m)}\psi^m_{i, l}\Big).
\]
Here $(\vx)_k$ is the $k^{th}$ component of any vector $\vx \in \R^D$.  We define $\balpha_k = \begin{bmatrix}(\va_1)_k & \cdots & (\va_n)_k \end{bmatrix}^\intercal \in \R^n$, with $A_k \in \R^{n \times n}$ and $\vb_k \in \R^{n}$ given as
\[
A_k(i, j) \coloneqq \frac{1}{2M}\sum_{l, m = 1}^{L - 1, M}\Big(\frac{\psi^m_{i, l}\psi^m_{j, l}}{\sigma_k^2(\rvx_l^m)}\Delta t_l\Big), \quad \vb_k(i) \coloneqq \frac{1}{2M}\sum_{l, m = 1}^{L - 1, M}\frac{\psi^m_{i, l}(\Delta\rvx_l^m)_k}{\sigma_k^2(\rvx_l^m)}.
\]
Then the definition in (\ref{eq:original_loss_dis_coeff}) can be rewritten as $\loss_{L, M, \gH}(\{\va_{\eta}\}_{i = 1}^n) = \sum_{k = 1}^D(\balpha_k^\intercal A_k\balpha_k - 2\balpha_k^\intercal\vb_k)$. Since each $\balpha_k^\intercal A_k\balpha_k - 2\balpha_k^\intercal\vb_k$ is decoupled from each other, we just need to solve simultaneously $A_k\hat\balpha_k - \vb_k = 0$, for $k = 1, \ldots, D$. Then we can obtain $\hat\vf(\vx) = \sum_{i = 1}^n\hat\va_i\psi_k(\vx)$. 
However when $\Sigma$ does not have a diagonal structure, we will have to resolve to gradient descent methods to minimize \eqref{eq:original_loss_dis_coeff} in order to find the coefficients $\{\va_i\}_{i = 1}^n$ for a total number of $nd$ parameters.  

If a data-driven basis is desired, we set $\gH$ to be the space of neural networks with the same depth, number of neurons, and activation functions in the hidden layers. Furthermore, we find $\hat\vf$ by minimizing the loss given by the definition in (\ref{eq:original_loss_dis}) using any deep learning optimizer, such as Stochastic Gradient Descent or Adam, from well-known deep learning packages.
\subsection{Proof of the Theorem}\label{append:proof}
We present the following definition about two different convergences of random variables.
\begin{definition}
A sequence $\{\rx_1, \rx_2, \cdots, \rx_n\}$ of scalar random variables, with cumulative distribution functions, $\{F_1, F_2, \cdots, F_n\}$, is said to converge in distribution to a random variable $\rx$ with cumulative distribution function $F$ if
\[
\lim_{n \rightarrow \infty} F_n(x) = F(x),
\]
for every number $x \in \R$ at which $F$ is continuous.  We denote such convergence as
\[
\rx_n \xrightarrow{D} \rx.
\]
We say $\rx_n$ convergences to $\rx$ in probability if for any $\epsilon > 0$, we have
\[
\lim_{n \rightarrow \infty} \mathbb{P}(|\rx_n - \rx| > \epsilon) = 0.
\]
We denote such convergence as
\[
\rx_n \xrightarrow{P} \rx.
\]
\end{definition}
The following lemma is needed for the convergence theorem.
\begin{lemma} \label{lemma:cont_inv}
Consider the space $(\mathbb{S}_{++}^n, \|\cdot\|_F)$ with $\mathbb{S}_{++}^n$ being the set of all $n \times n$ SPD matrices and $\norm{\cdot}_F$ denoting the Frobenius norm, then the inversion map $g: \mathbb{S}_{++}^n \to \mathbb{S}_{++}^n$ defined by $g(\mA) = \mA^{-1}$ for $\mA \in \mathbb{S}_{++}^n$ is continuous.
\end{lemma}
\begin{proof}
For any $\mA \in \mathbb{S}_{++}^n$ with $\det(\mA) > 0$, we have
\[
\mA^{-1} = \frac{\operatorname{adj}(\mA)}{\det(\mA)},
\]
where $\operatorname{adj}(\mA)$ is the adjugate matrix of $A$.  Each entry of $\operatorname{adj}(\mA)$ is a polynomial in the entries of $\mA$, and $\det(\mA)$ is also a polynomial in the entries of $\mA$. Since polynomials are continuous, both maps $\mA \mapsto \operatorname{adj}(\mA)$ and $\mA \mapsto \det(\mA)$ are continuous on $\R^{n \times n}$.  For $\mA \in \mathbb{S}_{++}^n$, we have $\det(\mA) > 0$, so the map $\mA \mapsto \frac{\operatorname{adj}(\mA)}{\det(\mA)}$ is continuous at $\mA$ as the composition of continuous functions.  Therefore, $g$ is continuous on $\mathbb{S}_{++}^n$.
\end{proof}
We present the following uniform law of large numbers theorem.  For the proof, please see~\citep{NEWEY19942111}.
\begin{theorem}[Uniform Law of Large Numbers~\citep{NEWEY19942111}]\label{thm:ulnn}
Let $\{\rx_i\}_{i=1}^\infty$ be i.i.d. and let $f(x, \theta)$ be some function defined for $\theta \in \Theta$. Assume:
\begin{enumerate}
\item $\Theta$ is compact;
\item for almost every $x$, the map $\theta\mapsto f(x,\theta)$ is continuous on $\Theta$, and for each $\theta\in\Theta$ the map $x\mapsto f(x,\theta)$ is measurable;
\item there exists a dominating function $h$ such that $\mathbb E[h(\rx)]<\infty$ such that $\|f(x,\theta)\|\le h(x)$ for all $\theta\in\Theta$.
\end{enumerate}
Then $\theta\mapsto \mathbb E[f(\rx,\theta)]$ is continuous in $\theta$ and
\[
\sup_{\theta\in\Theta}\norm{\frac{1}{n}\sum_{i=1}^n f(\rx_i,\theta)\;-\;\mathbb E[f(\rx,\theta)]} \xrightarrow{\,P\,} 0.
\]
\end{theorem}
The following theorem is needed to show convergence of vector-valued random variables.  For the proof, please see~\citep{Vaart_1998}.
\begin{theorem}[Theorem $5.9$ in~\citep{Vaart_1998}]\label{thm:varrt}
Let $\Psi_n:\Theta\to\mathbb R^k$ be random vector–valued functions and
$\Psi:\Theta\to\mathbb R^k$ a fixed vectored valued function of $\theta$.
Suppose that for every $\varepsilon>0$:
\[
\sup_{\theta\in\Theta}\,\|\Psi_n(\theta)-\Psi(\theta)\|\xrightarrow{\,P\,}0, \qquad \inf_{\theta:\, \norm{\theta - \theta_0}\ge \varepsilon}\|\Psi(\theta)\| \;>\; 0
\;=\; \|\Psi(\theta_0)\|.
\]
Then any sequence of estimator \(\hat \theta_n\) such that \(\Psi_n(\hat\theta_n) = o_p(1)\) converges in probability to \(\theta_0\). 
\end{theorem}

We are now ready to show the proof of the convergence theorem.
\begin{proof}
We need to introduce a few quantities before we can establish the proof.  First, we introduce the continuous form of $\loss_M$.  As $M \rightarrow \infty$, by law of large numbers, we have
\[
\lim_{M \rightarrow \infty}\loss_M(\tilde\vf) = \loss_{\infty}(\tilde\vf) = \frac{1}{2}\E\big[\int_0^T\inp{\tilde\vf_t, \,(\Sigma_t)^{-1}\tilde\vf_t} \, \dt - 2\int_0^T\inp{\tilde\vf_t, \,(\Sigma_t)^{-1}\dif\rvx_t}\big],
\]
where $\tilde\vf_t = \tilde\vf(\rvx_t)$, $\Sigma_t = \Sigma(\rvx_t)$.  When $\sH$ is finite dimensional, then for any $\tilde\vf \in \sH$, we have
\[
\tilde\vf(\vx) = \sum_{\eta = 1}^n\alpha_{\eta}\psivec_{\eta}(\vx) = \Psi(\vx)\alphavec, \quad \alphavec = \begin{bmatrix} \alpha_1 \\ \vdots \\ \alpha_n\end{bmatrix}.
\]
Therefore, the two losses can be re-written as
\[
\begin{aligned}
    \loss_M(\tilde\vf) &= \frac{1}{2M}\sum_{m = 1}^M\Big(\int_0^T(\Psi_t^m\alphavec)^\intercal(\Sigma_t^m)^{-1}\Psi_t^m\alphavec \, dt - 2\int_0^T(\Psi_t^m\alphavec)^\intercal(\Sigma_t^m)^{-1}\dif\rvx_t^m\Big), \\
    \loss_{\infty}(\tilde\vf) &= \frac{1}{2}\E\Big[\int_0^T(\Psi_t\alphavec)^\intercal(\Sigma_t)^{-1}\Psi_t\alphavec \, dt - 2\int_0^T(\Psi_t\alphavec)^\intercal(\Sigma_t)^{-1}\dif\rvx_t\Big), \\    
\end{aligned}
\]
Abusing the notation, we will use $\loss_M(\tilde\vf)$ and $\loss_M(\alphavec)$ interchangeably; similarly for $\loss_{\infty}(\tilde\vf)$ and $\loss_{\infty}(\alphavec)$, since $\alphavec$ and $\tilde\vf$ have a one-on-one correspondence once a $\sH$ is chosen.

\noindent Next, we will assume the following
\[
\begin{cases}
    &\E[\int_0^T\norm{\Psi_t^\intercal\Sigma_t^{-1}\Psi_t}_2 \, dt] < \infty, \\
    &\E[\int_0^T\norm{\Psi_t^\intercal\Sigma_t^{-1}\vf(\rvx_t)}_2 \, dt] < \infty, \\
    &\E[\int_0^T\norm{\Psi_t^\intercal\sigma_t^{-1}}_2 \, dt] < \infty, \\
\end{cases}
\]
Differentiating $\loss_M$ w.r.t to $\alphavec$ gives
\[
\nabla_{\alphavec}\loss_M(\alphavec) = \frac{1}{M}\sum_{m = 1}^M\big(\int_0^T(\Psi_t^m)^\intercal(\Sigma_t^m)^{-1}(\Psi_t^m\alphavec \, \dt - \dif\rvx_t^m)\big).
\]
Let
\[
\begin{aligned}
\phivec_m(\alphavec) &\coloneqq \int_0^T (\Psi_t^m)^\intercal(\Sigma_t^m)^{-1}(\Psi_t^m\alphavec\dt - \dif\rvx_t^m),\\
&= \int_0^T (\Psi_t^m)^\intercal(\Sigma_t^m)^{-1}(\tilde\vf_t^m\dt - \vf_t^m\dt - \sigma_t^m\dif\rvw_t^m), \\
&= \int_0^T (\Psi_t^m)^\intercal(\Sigma_t^m)^{-1}(\tilde\vf_t^m - \vf_t^m)\,\dt - \int_0^T(\Psi_t^m)^\intercal(\sigma_t^m)^{-1}\dif\rvw_t^m. 
\end{aligned}
\]
and define $\Phi_M(\alphavec) \coloneqq \frac{1}{M}\sum_{m = 1}^M\phivec_m(\alphavec)$.  First, by It\^{o}'s formula
\[
\E[\int_0^T(\Psi_t^m)^\intercal(\sigma_t^m)^{-1}\dif\rvw_t^m] = \vzero.
\]
Then
\[
\begin{aligned}
\E[\phivec_m(\alphavec)] &= \E[\int_0^T (\Psi_t^m)^\intercal(\Sigma_t^m)^{-1}(\tilde\vf_t^m - \vf_t^m)\,\dt], \\
&= \E[\int_0^T \Psi_t^\intercal\Sigma_t^{-1}(\tilde\vf_t - \vf_t)\,\dt]
\end{aligned}
\]
Define
\[
\Phi_{\infty}(\alphavec) = \lim_{m \rightarrow \infty}\Phi_M(\alphavec) = \E[\phivec_m(\alphavec)].
\]
By theorem \ref{thm:ulnn}, since $\sH$ is compact, $\phivec_m$ is continuous at each $\alphavec$ and it is also bounded (by one of our assumptions).  Moreover
\[
\sup_{\tilde\vf \in \sH}\norm{\Phi_M(\alphavec) - \Phi_{\infty}(\alphavec)} = \sup_{\tilde\vf \in \sH}\norm{\frac{1}{M}\sum_{m = 1}^M\phivec_m(\alphavec) - \E[\Phi_m(\alphavec)]} \xrightarrow{P} 0.
\]
Since $\vf \in \sH$, then $\vf(\vx) = \Psi(\vx)\alphavec_f$, then
\[
\begin{aligned}
\Phi_{\infty}(\alphavec) &= \E[\int_0^T \Psi_t^\intercal\Sigma_t^{-1}(\tilde\vf_t - \vf_t)\,\dt], \\
&= \E[\int_0^T \Psi_t^\intercal\Sigma_t^{-1}(\Psi_t\alphavec - \Psi_t\alphavec_f)\, \dt], \\
&= \E[\int_0^T \Psi_t^\intercal\Sigma_t^{-1}\Psi_t\, \dt](\alphavec - \alphavec_f)\, \\
&= \mA(\alphavec - \alphavec_f).
\end{aligned}
\]
Since $\mA$ is SPD, Let $\lambda_{\min}(\mA) > 0$ be the minimal eignevalue of $\mA$, then for all $\tilde\vf \in \sH$, 
\[
\norm{\Phi_{\infty}(\alphavec)} = \norm{\mA(\alphavec - \alphavec_f)} \ge \lambda_{\min}(\mA)\norm{\alphavec - \alphavec_f}.
\]
Therefore, for any $\epsilon > 0$, we have
\[
\inf_{\norm{\alphavec - \alphavec_f} \ge \epsilon}\norm{\Phi_{\infty}(\alphavec)} \ge \inf_{\norm{\alphavec - \alphavec_f} \ge \epsilon}\lambda_{\min}(\mA)\norm{\alphavec - \alphavec_f} \ge \lambda_{\min}(\mA)\epsilon > 0, 
\]
observe that $\Phi_{\infty}(\alphavec_f) = \vzero$.  By theorme \ref{thm:varrt}, we conclude that
\[
\hat\vf_M \xrightarrow{P} \vf, \quad \text{convergence in probability}. 
\]

\noindent Next, recall
\[
\Phi_M(\alphavec) = \frac{1}{M}\sum_{m = 1}^M\int_0^T(\Psi_t^m)^\intercal(\Sigma_t^m)^{-1}(\Psi_t^m\alphavec\dt - \dif\rvx_t^m),
\]
define
\[
\mA_M = \frac{1}{M}\sum_{m = 1}^M\int_0^T(\Psi_t^m)^\intercal(\Sigma_t^m)^{-1}\Psi_t^m \, dt.
\]
Since $\vf(\vx) = \Psi(\vx)\alphavec_f$, hence
\[
\begin{aligned}
    \phivec_m(\alphavec_f) &= \int_0^T(\Psi_t^m)^\intercal(\Sigma_t^m)^{-1}(\Psi_t^m\alphavec_f\dt - \dif\rvx_t^m), \\
    &= \int_0^T(\Psi_t^m)^\intercal(\Sigma_t^m)^{-1}(\vf_t^m\dt - \dif\rvx_t^m), \\
    &= -\int_0^T(\Psi_t^m)^\intercal(\Sigma_t^m)^{-1}\sigma_t^m\dif\rvw_t^m, \\
    &= -\int_0^T(\Psi_t^m)^\intercal(\sigma_t^m)^{-1}\dif\rvw_t^m
\end{aligned}
\]
This It\^{o} integral is square-integrable, and $\E[\phivec_m(\alphavec_f)] = \vzero$, and by It\^{o} isometry
\[
\text{Var}(\phivec_m(\alphavec_f)) = \E[\int_0^T\Psi_t^\intercal\Sigma_t^{-1}\Psi_t \, dt] = \mA < \infty.
\]
Sincr $\rvx_t^m$ is i.i.d, $\phivec_m(\alphavec_f)$ is also i.i.d.  Therefore, by the multivariate Central Limit Theorem, we have
\[
\sqrt{M}\Phi_M(\alphavec_f) = \frac{1}{\sqrt{M}}\sum_{m = 1}^M\phivec_m(\alphavec_f) \xrightarrow{D} \normal(\vzero, \mA).
\]
Furthermore, we also have the following (recall $\hat\vf(\vx) = \Psi(\vx)\hat\alphavec$)
\[
\Phi_M(\hat\alphavec) - \Phi_M(\alphavec_f) = \mA_M(\hat\alphavec - \alphavec_f), 
\]
Since $\Phi_M(\hat\alphavec) = \vzero$, we obtain
\[
\sqrt{M}(\hat\alphavec - \alphavec_f) = \sqrt{M}\mA_M^{-1}\Phi_M(\alphavec_f).
\]
Each entry of $\mA_M$ is square-integrable and by law of large numbers $\mA_M \rightarrow \mA$ as $M \rightarrow \infty$ in probability entrywise, hence
\[
\norm{\mA_m - \mA}_F \xrightarrow{P} 0.
\]
By lemma \ref{lemma:cont_inv}, the inversion mapping is continuous, hence
\[
\mA_M^{-1} \xrightarrow{P} \mA^{-1}.
\]
Putting them all together and by Slutsky's theorem, we end up with
\[
\sqrt{M}(\hat\alphavec - \alphavec_f) \xrightarrow{D} \normal(\vzero, \mA^{-1}).
\] 
Furthermore, for a fixed $\vx$, since $\hat\vf(\vx) = \Psi(\vx)\hat\alphavec$ and $\vf(\vx) = \Psi(\vx)\alphavec_f$, we finally have
\[
\sqrt{M}(\hat\vf - \vf) \xrightarrow{D} \normal(\vzero, \mA^{-1}). 
\] 
\end{proof}
\section{Examples}
In this section, we discuss the additional details for setting up the numerical examples and show additional examples.  In all examples, we use fairly complex covariance matrices, i.e., state-dependent matrices, in order to showcase the effectiveness of our learning.  The drift and noise estimations are carried out in both basis method and deep learning method with \ref{eq:original_loss} and \ref{eq:sigma_loss} being loss functions for estimating drift and covariance, respectively. The observations, serving as the input dataset for testing our method, are generated by the Euler-Maruyama scheme \cite{doi:10.1137/S0036144500378302}, utilizing the drift functions as we just mentioned. The basis space $\gH$ is constructed employing either B-spline or piecewise polynomial methods for maximum degree p-max equals $2$. For higher order dimensions where $d \geq 2$, each basis function is derived through a tensor grid product, utilizing one-dimensional basis defined by knots that segment the domain in each dimension. 

The parameters will be specified in each subsection of examples. The estimation results are evaluated using several different metrics. We record the noise terms, $\dif \rvw_t$, from the trajectory generation process and compare the trajectories produced by the estimated drift functions, $\hat{\vf}$, under identical noise conditions. We examine trajectory-wise errors using equation $\rho(\rvx) = \E\Big[\frac{1}{T}\int_{t = 0}^T\delta_{\rvx_t}(\rvx)\Big]$ with relative trajectory error and plot both $\vf$ and $\hat{\vf}$ to calculate the relative $L^2(\rho)$ error using \ref{eq:rho_norm}. When plotting, trajectories with different initial conditions are represented by distinct colors. In trajectory-wise comparisons, black solid lines depict the true trajectories, while blue dashed lines represent those generated by the estimated drift functions. Additionally, the empirical measure $\rho$ is shown in the background of each 1d plot. Furthermore, we assess the distribution-wise discrepancies between observed and estimated results, computing the Wasserstein distance at various time steps with \eqref{eq:wass_dis}. 
\subsection{Example: Interacting Particle Systems (IPS)}
 If we use the vectorized notations, i.e. 
\[
\rvx = \begin{bmatrix} \rvx_1 \\ \vdots \\ \rvx_N\end{bmatrix} \quad \text{and} \quad \rvw = \begin{bmatrix} \rvw_1 \\ \vdots \\ \rvw_N\end{bmatrix} \in \R^{D = Nd},
\]
and  
\[
\vf_{\phi}(\rvx) = \begin{bmatrix} \frac{1}{N}\sum_{j = 2}^N\phi(\norm{\rvx_j - \rvx_1})(\rvx_j - \rvx_1) \\ \vdots \\ \frac{1}{N}\sum_{j = 1}^{N - 1}\phi(\norm{\rvx_j - \rvx_N})(\rvx_j - \rvx_N)\end{bmatrix},  \quad
\sigma = \begin{bmatrix} \sigma^{\rx}(\rvx_1) & \zero & \cdots & \zero \\ \zero & \sigma^{\rx}(\rvx_2) & \cdots & \zero \\ \vdots & \vdots & \ddots & \vdots \\ \zero & \zero & \cdots & \sigma^{\rx}(\rvx_N) \end{bmatrix}. 
\]
Here each $\zero$ is a $d \times d$ matrix, $\vf: \R^D \rightarrow \R^D$ and $\tilde{\sigma}:\R^D \rightarrow \R^{D \times D}$.  Then the system can be put into one single SDE of the form $\dif \rvx_t = \vf(\rvx_t)\dif t + \tilde\sigma(\rvx_t)\dif \rvw_t$.  We will consider a weighted $\ell_2$ inner product for these vectors, i.e. for $\rvu, \rvv \in \R^d$ with
\[
\rvu = \begin{bmatrix} \rvu_1 \\ \vdots \\ \rvu_N\end{bmatrix}, \quad \rvv = \begin{bmatrix} \rvv_1 \\ \vdots \\ \rvv_N\end{bmatrix}, \quad \rvu_i, \rvv_i \in \R^{d}
\]
then
\[
\langle\rvu, \rvv\rangle_N = \frac{1}{N}\sum_{i = 1}^N\langle\rvu_i, \rvv_i\rangle, \quad \norm{\rvu}_N^2 = \langle\rvu, \rvu\rangle_N.
\]
With this new norm, we can carry out the learning as usual in $\R^d$ yet with a lower dimensional structure for $\vf_{\phi}$ and $\sigma^{\vx}$. With this setup, the loss of the noise in \eqref{eq:sigma_loss} will become
\[
\mathcal{E}_{\sigma}(\tilde \Sigma) = \E\Big[\frac{1}{N}\sum_{i = 1}^N\big([\rvx_{i, t}, \rvx_{i, t}]_0^T - \int_{t=0}^T (\tilde\sigma^\vx(\rvx_{i, t}))^2\dt\big)^2\Big],
\]
where we learn $\tilde\Sigma^\vx = (\tilde\sigma^\vx)^2$ as one single SPD matrix using the Cholesky decomposition method described in section \ref{sec:DL_learn}, and then take $\tilde\sigma^\vx = \sqrt{\tilde\Sigma^\vx}$. Next, the loss of the drift will become
\[
\loss_{\gH}(\varphi) = \frac{1}{2}\E\Big[\int_{t = 0}^T\langle\vf_\varphi(\rvx_t), \Sigma^{\dagger}(\rvx_t)\vf_\varphi(\rvx_t)\rangle_N \, \dt - 2\langle\vf_\varphi(\rvx_t), \Sigma^{\dagger}(\rvx_t)\dif \rvx_t\rangle_N\big)\Big].
\]
The two terms with the weighted $\ell_2$ inner product can be rewritten as
\[
\langle\vf_\varphi(\rvx_t), \Sigma^{\dagger}(\rvx_t)\vf_\varphi(\rvx_t)\rangle_N = \frac{1}{N^3}\sum_{i, j, k = 1}^N\varphi(r_{i, j, l}^m)\varphi(r_{i, k, l}^m)\langle\vr_{i, j, l}^m, (\tilde\sigma^\vx(\rvx_{i, l}^m))^{-2}\vr_{i, k, l}^m\rangle
\]
and
\[
\langle\vf_\varphi(\rvx_t), \Sigma^{\dagger}(\rvx_t)\dif \rvx_t\rangle_N = \frac{1}{N^2}\sum_{i, j = 1}^N\varphi(r_{i, j, l}^m)\langle\vr_{i, j, l}^m, (\tilde\sigma^\vx(\rvx_{i, l}^m))^{-2}(\rvx_{i, l + 1}^m - \rvx_{i, l}^m)\rangle,
\]
where $\rvx_{i, l}^m = \rvx_i^m(t_l)$, $\vr_{i, j, l}^m = \vx_j^m(t_l) - \vx_i^m(t_l)$, and $r_{i, j, l}^m = \norm{\vr_{i, j, l}^m}$.  In estimating $\phi$, we use \([r_{\min},r_{\max}]\) as the domain for estimation with \(r_{ij}=\|\rvx_j-\rvx_i\|\) represents the pairwise distance. We use a piecewise local B-spline basis of order up to \(3\) on domain \([r_{\min},r_{\max}]\).
Let \(\gH = \text{span}(\{\psi_\eta\}_{\eta=1}^{n})\) denote the associated compactly supported basis functions.
Then an estimator $\varphi \in \mathcal{H}$ has the form, i.e., $\varphi(r)=\sum_{\eta=1}^{n} a_\eta\,\psi_\eta(r)$.   For $\rvx=(\rvx_1^\top,\dots,\rvx_N^\top)^\top$, define interaction features indexed by the particle $i=1,\dots,N$, then each $\eta^{th}$ column of $\Phi$ is given as
\[
\big(\Phi_\eta(\vx)\big)_i \;=\; \frac{1}{N}\sum_{j = 1, j\neq i}^N\psi_\eta\!\big(r_{i, j}\big)\,(\sigma^{\vx}(\vx_i))^{-1}\vr_{i, j}\in\R^{d},
\]
where $\rv_{i, j} = \vx_j - \vx_i$ and $r_{i, j} = \norm{\vr_{i, j}}$. For $M$ trajectories $\{\rvx_l^m\}_{l, m = 1}^{L, M}$, set $\Delta \rvx^{m}_l=\rvx^{m}_{t_{l+1}}-\rvx^{m}_{t_l}$. Then the loss function reduces to
\[
\loss(a)=\tfrac12\,a^\top A\,a-b^\top a,
\]
where
\[
\begin{cases}
A_{ij} &\coloneqq \frac1{M}\sum_{m,l}
\langle \Phi_i\!\big(\rvx^{m}_{l}\big),\,
\,\Phi_j\!\rangle\,\Delta t,\\[0.5em]
b_i   &\coloneqq \frac1{M}\sum_{m,l}
\langle \Phi_i\!\big(\rvx^{m}_{l}\big),\,
\,\sigma^{-1}(\rvx_l^m)\Delta \rvx^{m}_l\rangle .
\end{cases}
\]
The estimator is the solution of the normal equations $A\,\hat a=b$, hence $\hat\phi(r)=\sum_{i=1}^{n}\hat a_i\,\psi_i(r)$.
\subsection{Example: SPDE estimation}
For any $N\in\mathbb{N}$, let $H^{N}=\spn\{h_1,\dots,h_N\}$ and $P^{N}\colon H\to H^{N}$ the projection operator. Then denote $\rvu^N = P^N \rvu = \sum_{k=1}^{N} \rvu_k(t) h_k(\rvx)$ as the Fourier  approximation of the solution $\rvu$ by the first $N$ eigenmodes $\rvu_k(t)=(\rvu(t), h_k)_H$. The projected solution $\rvu^N$ of \eqref{eq:SPDE_general} satisfies the following finite–dimensional SDE 
\begin{equation}\label{eq:true projected SPDE}
  \dif\rvu^{N}(t, \rvx)
  =
  P^{N}\!\bigl(\theta(\rvx)\,\Delta\rvu(t,\rvx)\bigr)\,\dif t
  +
  \sigma P^{N}\dif\rvw(t,\rvx).
\end{equation}

Since eigenmodes are coupled together in term $\theta(x)\Delta\rvu(t,\rvx)$, $P^N$ does not commute with $\theta(x)$, and to overcome this we consider a Galerkin type projection, i.e.
 \[
 \begin{aligned}
 &\tilde\rvu^N(t, \rvx) = \sum_{k = 1}^N\tilde\rvu_k(t)h_k(\rvx) &\approx \sum_{k = 1}^{\infty}\rvu_k(t)h_k(\rvx) = \rvu(t, \rvx), 
 \end{aligned}
 \]
and we have 
\begin{equation}\label{eq: Galerkin }
  \dif\tilde\rvu^{N}(t, \rvx)
  \;=\;
  P^{N}\!\bigl(\theta(\rvx)\,\Delta\tilde\rvu^N(t,\rvx)\bigr)\,\dif t
  +
  \sigma\,P^{N}\dif\rvw(t,\rvx),
\end{equation}
that we write in a matrix form,
\begin{equation}\label{eq:OU_tilde}
  \dif\tilde\rvU_N(t)
  \;=\;
  -G_{N}(\theta)\Lambda_{N}\,\tilde\rvU_N(t)\,\dif t
  + \sigma Q_{N}\,\dif\rvW_N(t),
\end{equation}
where
\begin{align}
\tilde\rvU_N &= \bigl(\tilde\rvu_{1},\dots,\tilde\rvu_{N}\bigr)^{\!\intercal}, \\
\dif\rvW_N &= \bigl(\dif\rvw_{1}(t),\dots,\dif\rvw_{N}(t)\bigr)^{\!\intercal}, \\
\Lambda_N &= \operatorname{diag}(\lambda_{1},\dots,\lambda_{N}), \\
Q_N &= \operatorname{diag}(q_{1},\dots,q_{N}),
\end{align}
and where the matrix
\(G_N(\theta)\in\R^{N\times N}\) has entries
\[
[G_{N}(\theta)]_{jk}=\inp{\theta(x)h_k, h_j},
\qquad 1\le j,k\le N.
\]
Choose any finite dimensional function space $\gH^\theta_n$ with basis 
\(\{\psi_{i}\}_{i=1}^{n}\), and approximate  $\theta(x)$ with respect to this basis, 
\[
\theta(x)\approx\sum_{i=1}^{n}a_{i}\,\psi_{i}(x),
\qquad \va=(a_{1},\dots,a_{n})^{\!\intercal}\in\R^{n}.
\]
For each \(i\) define the deterministic matrices
\begin{equation}\label{eq:B_matrices}
[B^{(i)}_N]_{jk}=\inp{\psi_i h_k, h_j},
\qquad 1\le j,k\le N,
\end{equation}
so that
\begin{equation}\label{eq:C_expand}
G_N(\theta)\;\approx\;\sum_{i=1}^{n}\va_{i}\,B^{(i)}_N.
\end{equation}
Let \(\Sigma=\sigma^{2}Q_NQ_N^{\intercal}\).
By our method in section ~\ref{sec:Learning Framework},
\begin{equation}\label{eq:loglik}
\loss(\va)
=\frac12\int_{0}^{T}
        \tilde{\rvU}^{\intercal}_N\Lambda_N G_N(\va)\Sigma^{-1}C_N(\va)\Lambda_N \tilde\rvU_N \dif t +  \int_{0}^{T}
        \tilde{\rvU}^{\intercal}_N\Lambda_N G_N(\va)\Sigma^{-1} \dif \tilde\rvU_N.
\end{equation}
Here, we abused the notation of $G_N(\va)$, since $\theta$ is approximated on $\gH^\theta_n$, it has a unique representation in terms of $\va$ once $\gH^\theta_n$ is chosen, we then define $G_N(\va) \coloneqq G_N(\theta = \sum_{i = 1}^n a_i\psi_i$).  With the expansion defined by (\ref{eq:C_expand}), we obtain a mass matrix $A$ and RHS vector $\vb$ having entries given as
\[
A_{mk} =\int_{0}^{T} \tilde{\rvU}^{\intercal}_N\Lambda_N B_N^{(m)}\Sigma^{-1}B_N^{(k)}\Lambda_N \tilde\rvU_N \dif t, \quad \vb_{m}=\int_{0}^{T}
       \tilde{\rvU}^{\intercal}_N\Lambda_N B_N^{(m)}\Sigma^{-1} \dif \tilde\rvU_N,
\]
for \(1\le m,k\le n\).  $A$ is apparently symmetric positive definite.  Next, the loss becomes
\[
\loss(\va) = \frac{1}{2}\va^\intercal A\va + \vb^\intercal\va,
\]
thus minimizing the loss is equivalent to solving solving the linear system \(\nabla\loss(\va)=0\), which gives the estimation coefficient as $\hat\va = -A^{-1}\vb$.
\end{document}